\documentclass[a4paper,12pt]{article}
\usepackage[cp1251]{inputenc}
\usepackage[english]{babel}
\usepackage[tbtags]{amsmath}
\usepackage{amsfonts,amssymb,mathrsfs,amscd}

\usepackage{wrapfig}
\usepackage[dvips]{graphicx}
\newtheorem{theorem}{Theorem}[section]
\newtheorem{lemma}{Lemma}[section]
\newtheorem{definition}{Definition}[section]

\date{}
 \numberwithin{equation}{section}

\begin{document}

\large
 \centerline{\bf Non-Haar MRA  on  local fields of positive characteristic}
  \centerline{\bf S.\,F.~Lukomskii, A.M.Vodolazov}
\noindent
 N.G.\ Chernyshevskii Saratov State University\\
 LukomskiiSF@info.sgu.ru\\
 vam21@yandex.ru

\begin{abstract}
  We propose a simple method to construct integral
  periodic mask and
  corresponding  scaling step functions that generate
  non-Haar orthogonal MRA on the local field $ F^{(s)}$ of positive characteristic $p$.
  To  construct this mask we use two new ideas. First, we  consider
  local field as vector space over the finite field $GF(p^s)$.
  Second, we  construct scaling function by arbitrary tree that has $p^s$ vertices. By
  fixed prime number $p$ there exist $p^{s(p^s-2)}$ such trees.\\
   Bibliography: 16 titles.
\end{abstract}
\noindent
 keywords: local fields, multiresolution
analysis, wavelet bases, trees.


\section{Introduction}\label{s1}
  First results on the wavelet analysis on local fields received Chinese mathem\-aticians
  Huikun Jiang, Dengfeng Li, and Ning Jin  in the
  article \cite{JLJ}. They  introduced the notion of MRA on local
  fields, for the fields  $ F^{(s)}$
  of positive characteristic $p$  proved some simple properties  and
  gave an algorithm for constructing wavelets for a known scaling
  function. Using these results they constructed MRA and
  corresponding  wavelets for the case when a scaling function is the chara\-cteristic
  function of unit ball $\cal D $. Such MRA is   called usually "Haar MRA" \   and corresponding
  wavelets -- "Haar wavelets".
  In \cite{LJ} wavelet frame on local field are constructed, a necessary
  condition and  sufficient conditions for wavelet frame on
  local fields are given too. Biswaranjan Behera and Qaiser Jahan \cite{BJ1}
  constructed the wavelet packets associated with  MRA on local
  fields of positive characteristic. In the article \cite{BJ2} the same
  authors proved that a function $\varphi\in L^2(\mathbb F^{(s)})$ is
  a scaling function for MRA in $L^2( F^{(s)})$ if and only if
  \begin{equation}\label{eq1.1}
    \sum_{k\in \mathbb N_0}|\hat\varphi(\xi+u(k))|^2=1\  for\  a.e. \ \xi \in
    {\cal D},
  \end{equation}
  \begin{equation}\label{eq1.2}
    \lim\limits_{j\to\infty}|\hat\varphi(\mathfrak p^j\xi)|=1 \ for\  a.e.\  \xi \in
    F^{(s)},
  \end{equation}
  and there exists an integral periodic function $m_0 \in L^2(\cal
  D)$ such that
  \begin{equation}\label{eq1.3}
  \hat\varphi(\xi)=m_0(\mathfrak p\xi)\hat\varphi(\mathfrak p\xi)\  for \ a.e.\  \xi
  \in F^{(s)}
  \end{equation}
    where $\{u(k) \}$ is the set of shifts, $\mathfrak p$ is a prime element.
    B.Behera and Q.Jahan \cite{BJ3} proved also  if the translates of the scaling
    functions of two multiresolution analyses are biorthogonal, then the associated
    wavelet families are also biortho\-gonal.
  So, to construct MRA on  a local field $\mathbb F^{(s)}$ we must
  construct an integral periodic mask $m_0$ with conditions
  (\ref{eq1.1}-\ref{eq1.3}). To solve this problem using  prime element methods developed in \cite{MT} is not
  simply.       Currently there are no effective methods for constructing such masks
  and scaling functions. In articles \cite{JLJ}-\cite{BJ3}   only Haar wavelets are
  obtained.

  In this paper, we propose a simple method to construct integral periodic masks and
  corresponding  scaling step functions that generate
  non-Haar ortho\-gonal MRA.
  To  construct this mask we use two new ideas. First, we  consider
  local field as vector space over the finite field $GF(p^s)$.
  Second, we  construct a scaling function by arbitrary tree that have $p^s$ nodes. For
  fixed prime number $p$ there exist $p^{s(p^s-2)}$ such trees.

   By $s=1$ the additive group $F^{(1)+}$ is a Vilenkin
   group. Issues of constructing of MRA and wavelets   on
   Vilenkin groups may be found in \cite{PF}-\cite{SLuk3}.

  The simplest example of a local field of characteristic zero is the field of p-adic
  numbers. Issues of constructing  MRA and wavelets
  on the field of p-adic numbers can be found in
  \cite{SVK}-\cite{AES}.

 The paper is organized as follows. We  consider local field
 $\mathbb F^{(s)}$ as a vector space over the finite field $GF(p^s)$.
 Therefore, in section 2, we recall some concepts and facts from
 the theory of finite fields and  define the local field $\mathbb
 F^{(s)}$  of positive characteristic $p$ as a set of infinite
 sequences $a=({\bf a}_j)$, where ${\bf a}_j\in GF(p^s)$.

In section 3 we prove that local field $\mathbb F^{(s)}$ is a
vector space over finite field $GF(p^s)$.

In section 4 we prove that the set $X$ of all characters of  local
field $\mathbb F^{(s)}$ also form  a vector space over finite
field $GF(p^s)$ with product as internal operation and powering as
external operation. We define Rademacher functions, find a general
view of characters, and prove a basic property of Rademacher
functions.

In section 5 we discuss the refinable equation and its mask.

  In  section 6 we consider refinable equation
  $$
   \hat\varphi(\chi)=m_0(\chi)\hat\varphi(\chi{\cal A}^{-1})
  $$
  with step mask $m_0$ and find a necessary and sufficient condition under which an integral periodic function
  $m_0$ is a mask of some refinement equation.

  In  section 7 we define $(N,M)$ elementary sets.
  We prove if $E\subset \mathbb F^{(s)}$ is $(N,M)$ elementary set and
  $|\hat\varphi(\chi)|={\bf 1}_E(\chi)$ on $X$ then the system of
 shifts $(\varphi(x\dot-h))_{h\in H_0}$ is an orthonormal system.

 In  section 8 we reduce the problem of construction of step
 refinable function to construction of some tree.
 We  consider some
 special class of refinable functions $\varphi(\chi)$ for which
 $|\hat\varphi(\chi)|$ is a characteristic function of a set. We
 introduce such concepts as "a set generated by a tree" \ and "a
 refinable step function generated by a tree" \ and prove, that
 every rooted tree containing $p^s$ nodes generates a refinable step function that generate an
 orthogonal MRA on local field $ F^{(s)}$. For $p=s=2$ we give an
 example of  a refinable  step function that generate non-Haar MRA.

 Using the results of  the
  article \cite{JLJ} we can construct now corresponding
  wavelets. This example shows that MRA on local field gives an effective method to
construct multidimensional step wavelets.\\

\section{Preliminaries}\label{s2}
We will consider two objects: Vilenkin groups and local fields.
Let $p$ be a prime number. Vilenkin group $(\mathfrak G,\dot+)$
consists of sequences
$$
a=(a_n)_{n\in\mathbb
Z}=(\dots,a_{n-1},a_{n},a_{n+1},\dots),\;a_j=\overline{0,p-1},
$$
in which only a finite number of terms with negative numbers are
nonzero. The operation $\dot+$ is defined as component wise
addition modulo $p$, i.e.
$$
a\dot+b=(a_n)\dot+(b_n)=((a_n+b_n){\rm mod}\,p)_{n\in\mathbb Z}.
$$
The topology in $\mathfrak G$ is determined by subgroups
$$
\mathfrak G_n=\{a\in \mathfrak
G:\;a=(\dots,0_{n-1,}a_{n-1},a_{n},a_{n+1},\dots)\}.
$$
The equality
$$
\rho(a,b)=\left\{\begin{array}{ll}
  \frac{1}{p^n}; & a_m\ne b_n,\;a_j=b_j\;
  \mbox{for}\;j<n \\
  0; & a_j=b_j\;
  \mbox{for}\;j\in \mathbb Z \\
\end{array}\right.
$$
is the non-Archimedian distance on $(\mathfrak G,\dot+)$. If $\mu$
is the Haar measure on $\mathfrak G$ then $\mu(\mathfrak
G_n\dot+g)=\mu\mathfrak G_n=\frac{1}{p^n}$, $n\in\mathbb Z$. The
dilation operator $\mathcal{A}$ is defined by the equation
$$
 \mathcal{A}(a)=(b_n)_{n\in\mathbb Z},\;\;b_n=a_{n+1}.
$$
It is evident that $\mathcal{A}\mathfrak G_n=\mathfrak G_{n-1}$
and $\int\limits_{\mathfrak
G}f(\mathcal{A}u)\,d\mu=\frac{1}{p}\int\limits_{\mathfrak
G}f(x)\,d\mu$.

By a local field we will mean a field $K$ which is locally
compact, non-discrete and totally disconnected. We will consider
local fields with positive characteristic only. By
Pontrjagin--Kovalsky theorem \cite{GGP-Sh} such field is
isomorphic to the set $K_L(z)$ of formally Loran series
\begin{equation}\label{eq2.1}
\sum_{n=N}^\infty a_nz^n
\end{equation}
with ${\bf a}_n\in GF(p^s)$ where  $s\in \mathbb N$ and $p$ is a
prime number. Local field of positive characteristic is denote
$F^{(s)}$.

 Let $GF(p)$ be a ring (field) of residue class on modulo $p$.
  The finite field $GF(p^s)$ consist of vectors
 ${\bf a}=(a^{(0)},a^{(1)},\dots,a^{(s-1)})$, where  $a^{(j)}\in
 GF(p)$.
  The addition operation $(\bf a)\dot+(\bf
 b)$ is defined coordinate-wise i.e.
 $$
 {\bf a}\dot+{\bf b}=(a^{(j)}+b^{(j)})\,{\rm
 mod}\,p)_{j=0}^{s-1}.
 $$
 To define a product ${\bf a}{\bf b}$ it is necessary to represent
 vectors ${\bf a}$ and ${\bf b}$ as polynomials
 $$
{\bf a}=\sum_{j=0}^{s-1}a^{(j)}t^j, {\bf
b}=\sum_{j=0}^{s-1}b^{(j)}t^j
 $$
 and multiply  these polynomials over the field  $GF(p)$. We
 obtain the polynomial
 $$
 Q(t)=\sum_{j=0}^{s-1}\sum_{k=0}^{s-1}a^{(j)}b^{(k)}t^{j+k}=\sum_{l=0}^{2s-2}t^l\sum_{k,j:\,k+j=l}a^{(j)}b^{(k)}
 $$
in which coefficients
$\beta_l=\sum\limits_{k,j:\,k+j=l}a^{(j)}b^{(k)}$ are calculating
in the field $GF(p)$. Then we take a prime  polynomial $p_s(t)$ of
degree $s$ and divide polynomial $Q(t)$ by $p_s(t)$ over the field
$GF(p)$
 $$
 Q(t)=p_s(t)q(t)+H(t).
 $$
 Coefficients $b_0,b_1,\dots, b_{s-1}$ of this rest ${ H}(t)$ are components of product ${\bf a}{\bf b}$.
 It is know that a prime polynomial $p_s(t)$ over the field
 $GF(p)$ exists but not only one. A prime polynomial $p_s(t)$ can
 be found by exhaustion.

We return to local fields. The sum and product of Loran series
\eqref{eq2.1} are defined in the standard way, i.e. if
  $$
 a=\sum_{j=k}^{\infty}{\bf a}_jt^j, \;b=\sum_{j=k}^{\infty}{\bf b}_jt^j
 $$
then
 \begin{equation}\label{eq2.2}
 a\dot+b=\sum_{j=k}^{\infty}({\bf a}_j\dot+{\bf b}_j)t^j,\;{\bf a}_j\dot+{\bf b}_j=((a_j^{(l)}+
 b_j^{(l)}){\rm mod}\, p)_{l=0}^{s-1}
 \end{equation}
\begin{equation}\label{eq2.3}
 ab=\sum_{l=2k}^{\infty}t^l\sum_{j,\nu:\,j+\nu=l}{\bf a}_j {\bf
 b}_\nu.
 \end{equation}
Topology in $F^{(s)}$ is given by neighborhood basis of zero
   $$
   F^{(s)}_n=\left\{a=\sum_{j=n}^\infty{\bf a}_jt^j:\;{\bf a}_j\in GF(p^s)\right\}.
   $$
   If
   $$
   a=\sum_{j=n}^\infty {\bf a}_jt^j,\; {\bf a}_n\ne 0
   $$
then we put by definition  $\|a\|=\frac{1}{p^{sn}}$. Consequently
   $$
   F^{(s)}_n=\left\{ x\in F^{(s)}:\;\|x\|\le \frac{1}{p^{sn}}\right\}.
   $$
By $F^{(s)+}$  denote the additive group of field  $F^{(s)}$.
Neighborhoods $F^{(s)}_n$ are compact subgroups of group
   $F^{(s)+}$. We will denote them as $F^{(s)+}_n$. The next properties are fulfill \\
   1)$\dots\subset F^{(s)+}_1\subset F^{(s)+}_0\subset
   F^{(s)+}_{-1}\subset\dots$\\
   2)$F^{(s)+}_n/ F^{(s)+}_{n+1}\cong GF(p^s)$ and $\sharp (F^{(s)+}_n/ F^{(s)+}_{n+1})=p^s$.

Therefore we will assume that a local field $F^{(s)}$ of positive
characteristic consists of infinite sequences
  $$
  a=(\dots ,{\bf 0}_{n-1},{\bf a}_n,{\bf a}_{n+1},\dots,),\;{\bf a}_j=(a_j^{(0)},a_j^{(1)},\dots,a_j^{(s-1)})\in GF(p^s)
  $$
in which only finite number of element ${\bf a}_j$ with negative
numbers are nonzero. The sum and product are defined as
  \begin{equation}\label{eq2.4}
 a\dot+b=({\bf a}_j\dot+ {\bf b}_j)_{i\in \mathbb Z},\ {\bf a}_j\dot+ {\bf b}_j=(a_j^{(\nu)}+
 b_j^{(\nu)}{\rm mod}\,p)_{\nu=0}^{s-1},
 \end{equation}
 \begin{equation}\label{eq2.5}
 ab= (\sum_{i,j:i+j=l}({\bf a}_i{\bf b}_j))_{l\in \mathbb Z}
 \end{equation}
 In this case
 $$
 \|a\|=\|(\dots,{\bf 0}_{n-1},{\bf a}_n,{\bf a}_{n+1},\dots)\|=\frac{1}{p^{sn}}
 \; \mbox{\rm if}\;  {\bf a}_n\ne {\bf 0},
 $$

 $$
   F^{(s)}_n=\{a=({\bf a}_j)_{j\in \mathbb Z}:\;{\bf a}_j\in GF(p^s);\; {\bf a}_j={\bf 0}\; \forall j<n
   \},
   $$
   $$
   \dots\subset F^{(s)}_1\subset F^{(s)}_0\subset
   F^{(s)}_{-1}\subset\dots,
   $$
   $F^{(s)}_n$ -- are compact subgroups in $F^{(s)+}$ and
   $\sharp(F^{(s)}_n/F^{(s)}_{n+1})=p^s$.

   It follows that $F^{(1)+}$ is are Vilenkin group. The converse is true also: in Vilenkin group
    $(\mathfrak G,\dot +)$ we can define product by \eqref{eq2.5}.
    With such operation $(\mathfrak G,\dot+,\cdot)$ will be a
    field.
    Since $F^{(1)+}$ is are Vilenkin group, it follows that\\
 1) $\int\limits_{{\mathfrak G_0}^\bot}(\chi,x)\,d\nu(\chi)={\bf
 1}_{{\mathfrak G_0}}(x)$,\
 2) $\int\limits_{\mathfrak G_0}(\chi,x)\,d\mu(x)={\bf 1}_{\mathfrak G_0^\bot}(\chi)$.\\
  3) $\int\limits_{\mathfrak G_n^\bot}(\chi,x)\,d\nu(\chi)=p^n{\bf
  1}_{\mathfrak G_n}(x)$,\
  4) $\int\limits_{\mathfrak G_n}(\chi,x)\,d\mu(x)=\frac{1}{p^n}{\bf
  1}_{\mathfrak G_n^\bot}(\chi)$\\
 where $\mathfrak G_n=F_n^{(1)+}$.

 From the definition of $F^{(s)}$ it follows that additive group
 $F^{(s)+}$ is also Vilenkin group $\mathfrak G$ and $F_n^{(s)+}=\mathfrak G_{ns}$.
\section{Locally field of positive  characteristic as vector space over a finite field}
 Let $(\mathfrak G,\dot+)$ be a Vilenkin group. We can define the multiplication operation
 on a number $\lambda\in GF(p)$ by the equation
$$
a\lambda=\underbrace{a\dot+a \dot+\dots \dot+a}_\lambda.
$$
Define the modulus of $\lambda$  as
$$
 |\lambda|=\left\{
           \begin{array}{ll}
 1,& \lambda \ne 0,\\
 0,& \lambda = 0,\\
  \end{array} \right.
$$
and the norm of $a\in\mathfrak G$ by the equation
 \begin{equation}\label{eq3.1}
\|a\|=p^{-n}
 \end{equation}
if
$$
 a=(\dots 0_{n-1}a_na_{n+1}\dots), n\in\mathbb Z, a_j\in
 Z_p, a_n\ne 0.
$$
Since $GF(p)$ is a field, it follows that $(\mathfrak G,\dot+)$ is
a vector space over the field  $ GF(p)$ and the equation
\eqref{eq3.1} defines  a norm in $(G,\dot+,\cdot\lambda)$.

Now we consider local field $F^{(s)}$ with positive characteristic
$p$. Its elements are infinite sequences
$$
a=(\dots , {\bf 0}_{n-1}, {\bf a}_n, {\bf a}_{n+1},\dots), {\bf
a}_j\in GF(p^s)
$$
where
$$
 {\bf a}_j=(a_j^{(0)},a_j^{(1)},\dots,a_j^{(s-1)}),\ a_j^{(\nu)}\in
 Z_p.
$$
Let $\lambda\in GF(p^s)$. By the definition
$\|a\|=\frac{1}{p^{sn}}$ if ${\bf a}_n\ne {\bf 0}$. Since
$$
 \lambda a=(\dots {\bf 0}_{-1},\lambda,{\bf 0}_1,\dots)\cdot(\dots
{\bf 0}_{n-1},{\bf a}_n,{\bf a}_{n+1},\dots)=
$$
$$
 (\lambda+{\bf 0}x+{\bf 0}x^2+\dots)({\bf a}_nx^n+{\bf a}_{n+1}x^{n+1}+\dots)=\lambda {\bf a}_nx^n+\lambda
 {\bf a}_{n+1}x^{n+1}+\dots=
 $$
 $$
 =(\dots {\bf 0}_{n-1},\lambda {\bf a}_n,\lambda {\bf a}_{n+1},\dots)
$$
it follows that the product  $\lambda {\bf a}$ is defined
coordinate wise. With such operations   $F^{(s)}$ is a vector
space. If we define the modulus $|\lambda|$ by the equation
$$
 |\lambda|=\left\{
           \begin{array}{ll}
 1,& \lambda \ne 0,\\
 0,& \lambda = 0.\\
  \end{array} \right.
$$
and norm $\|a\|$ by the equation
\begin{equation}\label{eq3.2}
 \|a\|=\frac{1}{p^{sn}},\;{\bf a}_n\ne 0
\end{equation}
then we can consider the field  $F^{(s)}$ as a vector normalized
space under the field $GF(p^s)$.

For brevity we denote $K:=F^{(s)},\  K_n:=F^{(s)}_n$. Take an
element  $g\in K_1\setminus K_2$ and fixed it. It is known
\cite{MT} that any element  $a\in K$ may be written in the form
$$
a=\sum_{n\in \mathbb Z} \lambda_ng^n,\;\lambda_n\in U,
$$
where $U$ is a fixed full set of coset representatives $K_2$ in
$K_1$. We can prove a more general statement.
\begin{theorem}
Let $(g_n)_{n\in \mathbb Z}$ be a fixed basic sequence in $K$,
i.e. $g_n\in K_n\setminus K_{n+1}$. Any element $a\in K$ may by
written as sum of the series
\begin{equation}\label{eq3.3}
a=\sum_{n\in \mathbb Z} \overline{\lambda}_ng_n, \
\overline{\lambda}_n\in GF(p^s).
\end{equation}
\end{theorem}
{\bf Proof.} Let $a\in K$. If $a=0$ then the equation
\eqref{eq3.3} is evident. Let $a\ne 0$. Then exists $n\in \mathbb
Z$ such that $a\in K_n^+\setminus K_{n+1}^+$. It means that
$$
a=(\dots {\bf 0}_{n-1},{\bf a}_{n},{\bf a}_{n+1},\dots), \ {\bf
a}_j\in GF(p^s),{\bf a}_n\ne {\bf 0}.
$$
Show that there exists $\overline{\lambda}_n \in GF(p^s)$ such
that
$$
a= \overline{\lambda}_ng_n \dot+\alpha_{n+1},\  \alpha_{n+1}\in
K_{n+1}.
$$
Indeed, since $g_n\in K_n\setminus K_{n+1}$ it follows that
$$
g_n=(\dots {\bf 0}_{n-1},{\bf g}_{n}^{(n)},{\bf
g}_{n+1}^{(n)},\dots), \ {\bf g}_n^{(n)}\ne {\bf 0}.
$$

 Since $GF(p^s)$ is a field, it follows there exists $\overline{\lambda}_n \in GF(p^s)$
 such that $\overline{\lambda}_n{\bf g}_n^{(n)}={\bf a}_n$. Therefore
 $$
 \overline{\lambda}_ng_n=(\dots
{\bf 0}_{n-1},\overline{\lambda}_n{\bf
g}_{n}^{(n)},\overline{\lambda}_n{\bf g}_{n+1}^{(n)}\dots)= (\dots
{\bf 0}_{n-1},{\bf a}_n,\tilde{\bf a}_{n+1}\dots).
  $$
  Consequently
 $$
 a\dot-\overline{\lambda}_ng_n=(\dots
 {\bf 0}_{n-1},{\bf 0}_n,{\bf a}_{n+1}-\tilde{\bf a}_{n+1}\dots)=\alpha_{n+1}\in
K_{n+1}^+,
$$
i.e. $a=\overline{\lambda}_ng_n \dot+\alpha_{n+1}$. Continuing
this process,
we obtain \eqref{eq3.3}.  $\square$ \\
{\bf Corollary.} If $g\in
K_1\setminus K_2$ then $g^n\in K_n \setminus K_{n+1}$. Therefore
we can take $g_n=g^n$ in the equation \eqref{eq3.3}.
\begin{definition}
 The operator
 $$
 \mathcal{A}:a=\sum_{n\in \mathbb Z}\overline{\lambda}_ng_n\longmapsto \sum_{n\in \mathbb
 Z}\overline{\lambda}_ng_{n-1}
 $$
 is called a delation operator.
\end{definition}
{\bf Remark 1.} If $g_n=g^n$  and  $a=\sum\limits_{n\in \mathbb
 Z}\overline{\lambda}_ng^n$ then $ag^{-1}=\sum\limits_{n\in \mathbb
 Z}\overline{\lambda}_ng^{n-1}$. So the dilation operation may be defined by
 equation $\mathcal{A}x=g^{-1}x$.\\
 {\bf Remark 2.}  Since additive group
 $F^{(s)+}$ is Vilenkin group $\mathfrak G$ with $F_n^{(s)+}=\mathfrak G_{ns}$
 it follows that $\mathcal{A}K_n=\mathcal{A}K_{n-1} $ and
 $\int\limits_{K^+}f(\mathcal{A}u)\,d\mu=\frac{1}{p^s}\int\limits_{K^+}f(x)\,d\mu$.

\section{Set of characters as vector space over a finite field}
Since $F^{(s)+}$ is a Vilenkin group it follows that the set of
characters is a locally compact zero-dimensional group with
product as group operation
$$
(\chi\varphi)(a)=\chi(a)\cdot\varphi(a).
$$
 Denote the set of characters as $X$. We want to find the explicit form of characters.
 Let us define
the character $r_n$ in the following way. If
$$
a=(\dots , {\bf 0}_{k-1}, {\bf a}_k, {\bf a}_{k+1},\dots), {\bf
a}_j\in GF(p^s)
$$
and
$$
 {\bf a}_j=(a_j^{(0)},a_j^{(1)},\dots,a_j^{(s-1)}),\ a_j^{(\nu)}\in
 GF(p)
$$
 then  $r_n(a)=e^{\frac{2\pi i}{p}a_k^{(l)}}$, where $n=ks+l$, $0\leq l<s$.
 \begin{lemma}
Any character $\chi\in X$ can be expressed uniquely as product
\begin{equation}\label{eq4.1}
  \chi=\prod_{n=-\infty}^{+\infty}
  r_n^{\alpha_n}\;\;(\alpha_n=\overline{0,p-1}),
   \end{equation}
 in which the number of factors with positive numbers are finite.
 \end{lemma}
  {\bf Proof}. Let
  $$
   x=(...,{\bf 0},{\bf x}_j,...,{\bf x}_k,{\bf x}_{k+1},...),\
   {\bf x}_k=(x_{ks+0},x_{ks+1},...,x_{ks+(s-1)})
  $$
  Since $F^{(s)+}$ is a Vilenkin group, it follows that functions  $r_{ks+l}(x)=e^{\frac{2\pi
  i}{p}x_{ks+l}}$,
  are Rademacher functions on $F^{(s)+}$. Therefore any character
  $\chi$ may by expressed in the form \eqref{eq4.1}. $\square$
  \begin{definition}
  Write the character $\chi$ as
  $$
  \chi =\prod_{k\in \mathbb Z}r_{ks+0}^{a_k^{(0)}}
  r_{ks+1}^{a_k^{(1)}}\dots  r_{ks+s-1}^{a_k^{(s-1)}}
  $$
  and denote
  $$
  {\bf r}_k^{{\bf a}_k}:=r_{ks+0}^{a_k^{(0)}}
  r_{ks+1}^{a_k^{(1)}}\dots r_{ks+s-1}^{a_k^{(s-1)}},
  $$
  where ${\bf a}_k=(a_k^{(0)},a_k^{(1)},...,a_k^{(s-1)})\in
  GF(p^s)$. The function ${\bf r}_k={\bf r}_k^{(1,0,\dots,0)}$ is
  called Rademacher function.
 \end{definition}
 \begin{definition}
 Assume by the definition
  $$
   ({\bf r}_k^{{\bf a}_k})^{{\bf b}_k}:={\bf r}_k^{{\bf a}_k{{\bf b}_k}},\ {\bf a}_k, {\bf b}_k\in
  GF(p^s).
   $$
   \end{definition}
   In this case
      $$
   {\bf r}_k^{{\bf a}_k}=({\bf r}_k^{(1,0,...,0)})^{{\bf a}_k}={\bf r}_k^{(a_k^{(0)},a_k^{(1)},...,a_k^{(s-1)})}=
   r_{ks+0}^{a_k^{(0)}}r_{ks+1}^{a_k^{(1)}}...r_{ks+s-1}^{a_k^{(s-1)}}.
   $$
Therefore we can write $\chi$ as the product
\begin{equation}\label{eq4.2}
     \chi =\prod_{k\in \mathbb Z}{\bf r}_k^{{\bf a}_k}.
    \end{equation}
\begin{definition}
   Define $\chi^{{\bf b}}$, ${\bf b}\in GF(p^s)$ as
    $$
    \chi^{\bf b}:=\prod\limits_{k\in\mathbb Z} ({\bf r}_k^{{\bf a}_k})^{\bf
    b_k}.
    $$
\end{definition}
\begin{lemma} Let ${\bf r}_k$ be a Rademacher function. Then
 $$
 {\bf r}_k^{{\bf u}\dot+{\bf v}}={\bf r}_k^{{\bf u}}\cdot{\bf r}_k^{\bf
 v},\;{\bf u},{\bf v}\in GF(p^s).
 $$
 \end{lemma}
  {\bf Proof.} Using the definition of Rademacher functions we have for $x=(x_k^{(l)})$
 $$
 ({\bf r}_k^{{\bf u}}{\bf r}_k^{\bf v},x) =({\bf r}_k^{{\bf
 u}},x)({\bf r}_k^{{\bf
 v}},x)=\prod\limits_{l=0}^{s-1}e^{\frac{2\pi
 i}{p}u_{ks+l}^{(l)}x_k^{(l)}}\cdot\prod\limits_{l=0}^{s-1}e^{\frac{2\pi
 i}{p}v_{ks+l}^{(l)}x_k^{(l)}}=
 $$
 $$
 =\prod\limits_{l=0}^{s-1}e^{\frac{2\pi
 i}{p}(u_{ks+l}^{(l)}+v_{ks+l}^{(l)})x_k^{(l)}}=({\bf r}_k^{{\bf
 u}\dot+{\bf v}},x).\;\;\square
 $$
   \begin{theorem}
 The set of characters of the field $F^{(s)}$ is a vector space\\
   $(X,\; *,\; \cdot^{GF(p^s)})$ under the finite field  $GF(p^s)$ with product as interior
   operation and powering as exterior operation.
   \end{theorem}
 {\bf Proof.} 1) Check $\chi^{{\bf u}\dot+{\bf v}}=\chi^{{\bf
 u}}\chi^{{\bf v}}$ for ${\bf u},{\bf v}\in GF(p^s) $. Let
 $$
  \chi^{{\bf u}}=\prod\limits_{k\in\mathbb Z}{\bf r}_k^{{\bf
  a}_k{\bf u}},\;\;\chi^{{\bf v}}=\prod\limits_{k\in\mathbb Z}{\bf
  r}_k^{{\bf a}_k{\bf v}}.
 $$
Using lemma 4.2. we obtain
$$
\chi^{{\bf u}}\chi^{{\bf v}}=\prod\limits_{k\in\mathbb Z}{\bf
r}_k^{{\bf a}_k{\bf u}}{\bf r}_k^{{\bf a}_k{\bf
v}}=\prod\limits_{k\in\mathbb Z}{\bf r}_k^{{\bf a}_k({\bf
u}\dot+{\bf v})}=\chi^{{\bf u}\dot+{\bf v}}.
$$
2) Check the equation $\chi_1^{{\bf u}}\chi_2^{{\bf
u}}=(\chi_1\chi_2)^{{\bf u}}$. Let
$$
\chi_1^{{\bf u}}=\prod\limits_{k\in\mathbb Z}{\bf r}_k^{{\bf
a}_k{\bf u}},\;\chi_2^{{\bf u}}=\prod\limits_{k\in\mathbb Z}{\bf
r}_k^{{\bf b}_k{\bf u}}.
$$
Using lemma 4.2 we have
$$
\chi_1^{{\bf u}}\chi_2^{{\bf u}}=\prod\limits_{k\in\mathbb Z}{\bf
r}_k^{{\bf a}_k{\bf u}}\prod\limits_{k\in\mathbb Z}{\bf r}_k^{{\bf
b}_k{\bf u}}=\prod\limits_{k\in\mathbb Z}{\bf r}_k^{({\bf
a}_k\dot+{\bf b}_k){\bf u}}=(\chi_1\chi_2)^{{\bf u}}.
$$
3) Since the vector ${\bf 1}=(1,0,\dots,0)$ is a unity element of
multiplicative group of the field $GF(p^s)$ it follows that
$\chi^{{\bf 1}}=\chi^{(1,0,\dots,0)}=\prod\limits_{k\in\mathbb
Z}{\bf r}_k^{{\bf a}_k\cdot{\bf 1}}=\prod\limits_{k\in\mathbb
Z}{\bf
r}_k^{{\bf a}_k}=\chi$.\\
4)  The equality  $(\chi^{{\bf u}})^{{\bf v}}=\chi^{{\bf u}{\bf
v}}$ is true by the definition.

So, all axioms for exterior operation are fulfil. By lemma 4.2 all
axioms for interior operation are fulfil too. $\square$

   It follow from \eqref{eq4.2} that annihilator  $(F_k^{(s)})^\bot$ consist
   from characters of form $\chi={\bf r}_{k-1}^{{\bf a}_{k-1}}{\bf r}_{k-2}^{{\bf
   a}_{k-2}}...$. It is evident also that\\
   1) Rademacher system $({\bf r}_k)$ forms a basis of $(X,\; *,\; \cdot^{GF(p^s)})$,\\
   2) any sequences of characters $\chi_k\in (F_{k+1}^{(s)})^\bot\setminus(F_k^{(s)})^\bot$
   forms a basis of $(X,\; *,\; \cdot^{GF(p^s)})$.\\
   3)$(F_k^{(s)})^\bot=\bigsqcup\limits_{{\bf a}_{k-1}\in GF(p^s)}(F_{k-1}^{(s)})^\bot {\bf r}_{k-1}^{{\bf
   a}_{k-1}}$.\\
  The next lemma is the basic property  of Rademacher
  functions    on local field with positive characteristic.
   \begin{lemma} Let
   $g_j=(\dots,{\bf 0}_{j-1},(1,0,\dots,0)_j,{\bf
   0}_{j+1},\dots)\in F^{(s)}$, ${\bf a}_k,{\bf u}\in GF(p^s)$. Then
   $({\bf r}_k^{{\bf a}_k},{\bf u}g_j)=1$ for any $k\ne j$.
   \end{lemma}
   {\bf Proof.} Since ${\bf u}g_j=(\dots,{\bf 0}_{j-1},(u^{(0)},
   u^{(1)},\dots,u^{(s-1)})_j,{\bf 0}_{j+1},\dots)$, it follows
   that
   $$
   ({\bf r}_k^{{\bf a}_k},{\bf
   u}g_j)=\prod\limits_{l=0}^{s-1}e^{\frac{2\pi
   i}{p}a_k^{(l)}u^{(l)}}=\prod\limits_{l=0}^{s-1}e^{0}=1.\;\;\square
   $$

   \begin{definition}
   Define a dilation operator $\mathcal{A}$ on the set of characters by
    the equation $(\chi \mathcal{A},x)=(\chi, \mathcal{A}x)$.
   \end{definition}
{\bf Remark.}  Since additive group
 $F^{(s)+}$ is Vilenkin group,
 it follows that $g_j\mathcal{A}=g_{j+1}, \
    (K_n^+)^\bot\mathcal{A}=(K_{n+1}^+)^\bot $ and
 $\int\limits_{X}f(\chi\mathcal{A})\,d\nu=\frac{1}{p^s}\int\limits_Xf(\chi)\,d\nu$.

 \section{MRA on local fields of positive characteristic}
  We will use Rademacher function to construct MRA on local fields of positive
  characteristic. We will assume
  $$
  g_n=(\dots,{\bf 0}_{n-1},(1,0,\dots,0)_n,{\bf
   0}_{n+1},\dots).
   $$
  \begin{lemma}
  Let  $K=F^{(s)}$ be a local field with characteristic $p$. Then for any $n\in \mathbb Z$ \\
  1) $\int\limits_{(K_n^+)^\bot}(\chi,x)\,d\nu(\chi)=p^{sn}{\bf 1}_{K_n^+}(x)$,\\
  2) $\int\limits_{K_n^+}(\chi,x)\,d\mu(x)=\frac{1}{p^{sn}}{\bf 1}_{(K_n^+)^\bot}(\chi)$.
  \end{lemma}
 {\bf Proof.}
  First we prove the equation 1).
  Since $K^+$ is a zero-dimensional group, it follows
  $$
    \int\limits_{(K_0^+)\bot}(\chi,x)d\nu(\chi)={\bf
    1}_{K_0^+}(x),\quad
    \int\limits_{K_0^+}(\chi,x)d\mu(x)={\bf
    1}_{(K_0^+)^\bot}(\chi).
  $$
  By the definition of dilation operator
  $$
   \int\limits_Xf(\chi{\cal
   A})\,d\nu(\chi)=p^s\int\limits_Xf(\chi)\,d\nu(\chi),\;\;\;{\bf
   1}_{K_n^+}(x)={\bf 1}_{K^+_0}({\cal A}^nx).
  $$
   Using these equations we have
 $$
   \int\limits_{(K_n^+)^\bot}(\chi,x)\,d\nu(\chi)=\int\limits_X{\bf
   1}_{(K_n^+)^\bot}(\chi)(\chi,x)\,d\nu(\chi)=
 $$
 $$
   =p^{sn}\int\limits_X(\chi{\cal
   A}^n,x){\bf 1}_{(K_n^+)^\bot}(\chi{\cal A}^n)\,d\nu(\chi)=
 $$
 $$
  =p^{sn}\int\limits_X(\chi,{\cal A}^nx){\bf
  1}_{(K_0^+)^\bot}(\chi)\,d\nu(\chi)=p^{sn}{\bf 1}_{K^+_0}({\cal
  A}^nx)=p^{sn}{\bf 1}_{K^+_n}(x).
 $$
 The second equation is proved by analogy.  $\square$

 \begin{lemma}
 Let  $\chi_{n,l}={\bf r}_n^{{\bf a}_n}{\bf r}_{n+1}^{{\bf
 a}_{n+1}}\dots {\bf r}_{n+l}^{{\bf a}_{n+l}}$ be a character does
 not belong to $({K_n^+})^\bot$. Then
 $$
 \int\limits_{({K_n^+})^\bot\chi_{n,l}}(\chi,x)\,d\nu(\chi)=p^{ns}(\chi_{n,l},x){\bf
 1}_{K_n^+}(x).
 $$
\end{lemma}
{\bf Proof.} Denote $\mathfrak G_n :={K_n^+}$. By analogy with
previously we have
$$
\int\limits_{\mathfrak G_n^\bot\chi_{n,l}}(\chi,x)\,d\nu(\chi)=
\int\limits_X{\bf 1}_{\mathfrak
G_n^\bot\chi_{n,l}}(\chi)(\chi,x)\,d\nu(\chi)= \int\limits_X{\bf
1}_{\mathfrak G_n^\bot}(\chi)(\chi_{n,l}\chi,x)\,d\nu(\chi)=
$$
$$
\int\limits_{\mathfrak
G_n^\bot}(\chi_{n,l},x)(\chi,x)\,d\nu(\chi)=p^{ns}(\chi_{n,l},x)
{\bf 1}_{\mathfrak G_n}(x). \;\;\square
$$
 \begin{lemma}
 Let $h_{n,l}={\bf
 a}_{n-1}g_{n-1}\dot+{\bf a}_{n-2}g_{n-2}\dot+\dots\dot+{\bf a}_{n-l}g_{n-l}\notin
 K_n^+$. Then
 $$
 \int\limits_{K_n^+\dot+h_{n,l}}(\chi,x)\,d\mu(x)=\frac{1}{p^{ns}}(\chi,h_{n,l}){\bf
 1}_{({K_n^+})^\bot}(\chi).
 $$
 \end{lemma}
 This lemma is dual to lemma 5.2.

  \begin{definition}
 Let $M,N\in\mathbb N$.
 Denote by  ${\mathfrak D}_M(K_{-N})$ the set of step-functions
 $f\in L_2(K)$ such that 1)${\rm supp}\,f\subset K_{-N}$, and 2)
 $f$ is constant on cosets $K_M\dot+g$. Similarly is defined ${\mathfrak
 D}_{-N}(K_{M}^\bot)$.
 \end{definition}
\begin{lemma}
 Let $M,N\in\mathbb N$. $f\in \mathfrak D_M(K_{-N})$ if and only if $\hat f\in \mathfrak
 D_{-N}(K_M^\bot)$.
\end{lemma}
{\bf Proof.} It is evident since additive group $F^+$ is Vilenkin
group.

\begin{lemma}
Let $\varphi\in L_2(K)$. The system $(\varphi(x\dot-h))_{h\in
H_0}$ is orthonormal if and only if the system
$\left(p^{\frac{ns}{2}}\varphi({\cal A}^nx\dot-h)\right)_{h\in
H_0}$ is orthonormal.
\end{lemma}
{\bf Proof}. This lemma follows from the equation
$$
\int\limits_Kp^{\frac{ns}{2}}\varphi({\cal
A}^nx\dot-h)p^{\frac{ns}{2}}\overline{\varphi({\cal
A}^nx\dot-g)}\,d\mu=\int\limits_K\varphi(x\dot-h)\overline{\varphi(x\dot-g)}\,d\mu.\;\;\square
$$

\begin{definition}
  A family of closed subspaces $V_n$, $n\in\mathbb Z$,
 is said to be a~multi\-resolution analysis  of~$L_2(K)$
 if the following axioms are satisfied:
 \begin{itemize}
 \item[A1)] $V_n\subset V_{n+1}$;
 \item[A2)] ${\vrule width0pt
 depth0pt height11pt} \overline{\bigcup_{n\in\mathbb
 Z}V_n}=L_2(K)$ and $\bigcap_{n\in\mathbb Z}V_n=\{0\}$;
 \item[A3)] $f(x)\in V_n$  $\Longleftrightarrow$ \ $f({\cal A} x)\in V_{n+1}$ (${\cal A}$~is a~dilation
 operator);
 \item[A4)] $f(x)\in V_0$ \ $\Longrightarrow$ \
 $f(x\dot - h)\in V_0$ for all $h\in H_0$; ($H_0$ is analog of $\mathbb
 Z$).
  \item[A5)] there exists
 a~function $\varphi\in L_2(K)$ such that the system
 $(\varphi(x\dot - h))_{h\in H_0}$ is an orthonormal basis
 for~$V_0$.
\end{itemize}
 A function~$\varphi$ occurring in axiom~A5 is called
a~\textit{scaling function}.
\end{definition}
 It is clear that the axiom A5 follows the axiom A4.
Next we will follow the conventional approach. Let
 $\varphi(x)\,{\in}\, L_2(K)$, and suppose that
 $(\varphi(x\dot -\nobreak h))_{h\in H_0}$ is an~orthonormal
 system in~$L_2(K)$. With the function~$\varphi$ and the
 dilation operator~${\cal A}$, we define the linear subspaces
 $L_n=(\varphi({\cal A}x\dot - h))_{h\in H_0}$ and
 closed subspaces $V_n=\overline{L_n}$. It is evident that the functions
  $p^{\frac{ns}{2}}\varphi({\cal A}x \dot-h)_{h\in H_0}$ form
 an orthonormal basis for $V_n$, $n\in \mathbb Z$. Therefore the axiom A3 is fulfilled.
 If subspaces $V_j$ form
 a~MRA, then the function~$\varphi$ is said to \textit{generate}
 an~MRA in~$L_2(K)$. If a function $\varphi$ generates an MRA, then we obtain
 from the axiom A1
\begin{equation}
                                                                   \label{eq5.1}
  \varphi(x)=\sum_{h\in H_0}\beta_h\varphi({\cal
  A}x\dot-h)\;\;\left(\sum|\beta_h|^2<+\infty\right).
 \end{equation}
  Therefore we will look up a~function
 $\varphi\in L_2(K)$, which generates an~MRA
 in~$L_2(K)$, as a~solution of the refinement
 equation (\ref{eq5.1}). A solution of refinement equation (\ref{eq5.1}) is called a {\it refinable function}.
 \begin{lemma}
Let $\varphi \in \mathfrak D_M(K_{-N})$ be a solution of
(\ref{eq5.1}). Then
   \begin{equation}                                                 \label{eq5.2}
\varphi(x)=\sum_{h\in H_0^{(N+1)}}\beta_h\varphi({\cal A}x\dot-h)
 \end{equation}
\end{lemma}
 The proof repeats the proof of Lemma 4.1 in \cite{SLuk1}.

 \begin{theorem}
 Let $\varphi\in \mathfrak D_M(K_{-N})$ and let
 $(\varphi(x\dot-h))_{h\in H_0}$ be an orthonormal system.
 $V_n\subset V_{n+1}$ if and only if  the function $\varphi(x)$ is
 a solution of refinement equation (\ref{eq5.2}).
  \end{theorem}
  The proof repeats the proof of Theorem 4.2 in \cite{SLuk1}.
 \begin{theorem}{(\rm \cite{BJ2}, th.4.1).}
  Let $\varphi\in \mathfrak D_M(K_{-N})$ be a solution of the
  equation  (\ref{eq5.2}), $(\varphi(x\dot-h))_{h\in H_0}$ an
 orthonormal basis  in
 $V_0$
 Then  $\bigcap\limits_{n\in \mathbb Z}V_n=\{0\}$.
 \end{theorem}

 \begin{theorem}{\rm (\cite{BJ2}, th.4.3)}
  Let $\varphi\in \mathfrak D_M(K_{-N})$ be a solution of the  equation  (\ref{eq5.2}), $(\varphi(x\dot-h))_{h\in H_0}$ an orthonormal basis  in
 $V_0$, and $\hat\varphi (0)\neq 0$.
  Then \\ $\overline{\bigcup\limits_{n\in\mathbb Z}V_n}=L_2(K)$.
 \end{theorem}

 The refinement equation (\ref{eq5.2}) may be written in the form
 \begin{equation}                                                        \label{eq5.3}
 \hat\varphi(\chi)=m_0(\chi)\hat\varphi(\chi{\cal
  A}^{-1}),
 \end{equation}
  where

 \begin{equation}                                                        \label{eq5.4}
 m_0(\chi)=\frac{1}{p^s}\sum_{h\in
 H_0^{(N+1)}}\beta_h\overline{(\chi{\cal A}^{-1},h)}
 \end{equation}
 is a mask of the equation (\ref{eq5.3}).
 \begin{lemma}
 Let $\varphi\in\mathfrak D_M(K_{-N})$. Then the mask $m_0(\chi)$
 is constant on cosets $K_{-N}^\bot\zeta$.
 \end{lemma}
 {\bf Proof.} We will prove that $(\chi,{\cal A}^{-1}h)$ are
 constant on cosets $K_{-N}^\bot\zeta$. Without loss of generality,
 we can assume that $\zeta={\bf r}_{-N}^{{\bf a}_{-N}}\dots {\bf
 r}_{-N+s}^{{\bf a}_{-N+s}}\notin K_{-N}^\bot$. If
 $$
 h={\bf a}_{-1}g_{-1}\dot+\dots\dot+{\bf a}_{-N-1}g_{-N-1}\in
 H_0^{(N+1)},\ {\bf a}_j  \in GF(p^s)
 $$
 then
 $$
  {\cal A}^{-1}h={\bf a}_{-1}g_{0}\dot+\dots\dot+{\bf
  a}_{-N-1}g_{-N}\in
   K_{-N}.
 $$
 If $\chi\in K_{-N}^\bot\zeta$ then
 $\chi=\chi_{-N}\zeta$ where  $\chi_{-N}\in K_{-N}^\bot$. Therefore
 $(\chi,{\cal A}^{-1}h)=(\chi_{-N}\zeta,{\cal
   A}^{-1}h)=(\zeta,{\cal A}^{-1}h)$. This means that $(\chi,{\cal
   A}^{-1}h)$ depends on  $\zeta$ only. $\square$

 \begin{lemma}
 The mask  $m_0(\chi)$ is a periodic  function with any period
 ${\bf r}_1^{{\bf a}_1}{\bf r}_2^{{\bf a}_2}\dots {\bf r}_l^{{\bf a}_l}$ $(l\in\mathbb
 N,\; {\bf a}_j\in GF(p^s),\;j=\overline{1,l})$.
 \end{lemma}
 {\bf Proof.}
 Using the equation $({\bf r}_k^{{\bf b}_k},{\bf u}g_j)=1, (k\neq j)$ we find
 $$
 (\chi  {\bf r}_1^{{\bf b}_1}{\bf r}_2^{{\bf b}_2}\dots {\bf
 r}_l^{{\bf b}_l}  ,{\cal A}^{-1}h)= (\chi {\bf r}_1^{{\bf
 b}_1}{\bf r}_2^{{\bf b}_2}\dots {\bf r}_l^{{\bf b}_l}  ,
 {\bf a}_{-1}g_0\dot+{\bf a}_{-2}g_{-1}\dot+\dots\dot+{\bf a}_{-N-1}g_{-N})=
 $$
 $$
  =(\chi,{\bf a}_{-1}g_0\dot+{\bf a}_{-2}g_{-1}\dot+\dots\dot+{\bf a}_{-N-1}g_{-N})=(\chi{\cal
  A}^{-1},h).
 $$
 Therefore
 $m_0(\chi {\bf r}_1^{{\bf b}_1}{\bf r}_2^{{\bf b}_2}\dots {\bf
 r}_l^{{\bf b}_l})=m_0(\chi)$
  and the lemma  is proved. $\square$
 \begin{lemma}
  The mask $m_0(\chi)$ is defined by its values on cosets
  $K_{-N}^\bot {\bf r}_{-N}^{{\bf a}_{-N}}\dots {\bf r}_0^{{\bf a}_0}$
  $({\bf a}_j=(a_j^{(0)},a_j^{(1)},...,a_j^{(s-1)})\in GF(p^s) )$.
  \end{lemma}
 {\bf Proof.} Let us denote
 $$
  k=\sum_{j=0}^N(a_{-j}^{(0)}+a_{-j}^{(1)}p+...+a_{-j}^{(s-1)}p^{s-1})p^{sj}\in [0,p^{s(N+1)}-1],
 $$
 $$
   l=\sum_{j=1}^{N+1}(\alpha_{-j}^{(0)}+\alpha_{-j}^{(1)}p+\dots+
   \alpha_{-j}^{(s-1)}p^{s-1})p^{s(j-1)}\in [0,p^{s(N+1)}-1].
 $$
 Then (\ref{eq5.4}) can be written as the system
 \begin{equation}                                                            \label{eq5.5}
 m_0(\chi_k)=\frac{1}{p^s}\sum_{l=0}^{p^{s(N+1)}-1}\beta_l\overline{(\chi_k,{\cal
 A}^{-1}h_l)},\;k=\overline{0,p^{s(N+1)}-1}
  \end{equation}
 in the unknowns $\beta_l$. We consider the characters $\chi_k$ on
 the subgroup $K_{-N}^+$. Since ${\cal A}^{-1}h_l \in  K^+_{-N}$,
 it follows that the matrix
 $p^{-\frac{s(N+1)}{2}}\overline{(\chi_k,{\cal A}^{-1}h_l)}$ is
 unitary, and so the system (\ref{eq5.5}) has a unique solution for
 each finite sequence \\ $(m_0(\chi_k))_{k=0}^{p^{s(N+1)}-1}$.
 $\square$

 {\bf Remark.} The function $m_0(\chi)$
 constructing in  Lemma 5.9 may be not a mask for $\varphi\in
 \mathfrak D_M(\mathfrak G_{-N})$. In the section 4 we find
 conditions under which the function $m_0(\chi)$ will be a mask.
 \begin{lemma}
 Let $\hat f_0(\chi)\in{\mathfrak D}_{-N}(\mathfrak (K_1^+)^\bot)$.
 Then
  \begin{equation}                                                                       \label{eq5.6}
 \hat f_0(\chi)=\frac{1}{p^s}\sum_{h\in
 H_0^{(N+1)}}\beta_h\overline{(\chi,{\cal A}^{-1}h)}.
 \end{equation}
 \end{lemma}
 {\bf Prof.} Since
 $\int\limits_{(K_0^+)^\bot}(\chi,g)\overline{(\chi,h)}\,d\nu(\chi)=\delta_{h,g}$
 for $h,g\in H_0$ it follows that\\
 $\int\limits_{(K_1^+)^\bot}(\chi{\cal
 A}^{-1},g)\overline{(\chi{\cal
 A}^{-1},h)}\,d\nu(\chi)=p\delta_{h,g}$.
 Therefore we can
 consider the set $\left(\frac{{\cal
 A}^{-1}h}{\sqrt{p^s}}\right)_{h\in H_0^{(N+1)}}$ as an orthonormal
 system on $(K_1^+)^\bot$. We know (lemma 5.7) that $(\chi,{\cal
 A}^{-1}h)$ is a constant on cosets $(K_{-N}^+)^\bot\zeta$. It is
 evident the dimensional of $\mathfrak D_{-N}(\mathfrak
 (K_1^+)^\bot)$ is equal to $p^{s(N+1)}$. Therefore the system
 $\left(\frac{{\cal A}^{-1}h}{\sqrt{p^s}}\right)_{h\in H_0^{(N+1)}}$
 is an orthonormal basis for $\mathfrak D_{-N}((K_1^+)^\bot)$ and
 the equation (\ref{eq5.6}) is valid. $\square$

 \section{Non Haar wavelets}
 In this section we find the necessary and sufficient condition
 under which a step function $\varphi(x)\in {\mathfrak
 D}_M(K_{-N})$ generates an orthogonal MRA on the local field with positive characteristic.
   We  will prove also  that for any $n\in \mathbb N$ there exists a step function $\varphi$
  such that 1) $\varphi$ generate   an orthogonal MRA, 2) ${\rm supp}\,\hat{\varphi}\subset
 K_n^{\perp}$, 3) $\hat\varphi(K_{n}^\bot
  \setminus K_{n-1}^\bot)\not\equiv 0$.\\
 Note that the results of Sections 6, 7 and 8, there are analogues
 of the corresponding results for Vilenkin groups \cite{SLuk3}. Moreover, we use
 the same methods. This is possible since the basic property of
 Rademacher functions (Lemma 4.3) is satisfied.

 First we give a test under which the
 system of
  shifts $(\varphi(x\dot-h))_{h\in H_0}$ is an orthonormal system.
 \begin{theorem}
 Let  $\varphi(x)\in {\mathfrak
 D}_M(K_{-N})$. A shift's system
 $(\varphi(x\dot-h))_{h\in H_0}$ will be orthonormal if and only
 if for any
 $\overline{\alpha}_{-N},\overline{\alpha}_{-N+1},\dots,\overline{\alpha}_{-1}\in GF(p^s)$
 \begin{equation}                                                                \label{eq6.1}
 \sum_{\overline{\alpha}_{0},\overline{\alpha}_1,\dots,\overline{\alpha}_{M-1}\in GF(p^s)}
 |\hat\varphi((K_{-N}^+)^\bot
 {\bf r}_{-N}^{\overline{\alpha}_{-N}}\dots {\bf r}_0^{\overline{\alpha}_0}\dots
 {\bf r}_{M-1}^{\overline{\alpha}_{M-1}})|^2=1.
 \end{equation}
 \end{theorem}
 {\bf Proof.} First we prove that the system
 $(\varphi(x\dot-h))_{h\in H_0}$ will be orthonormal if and only
 if
 \begin{equation}                                                                   \label{eq6.2}
 \sum_{\overline{\alpha}_{-N},\dots,\overline{\alpha}_0,\dots,\overline{\alpha}_{M-1}}
 |\hat\varphi((K_{-N}^+)^\bot
 {\bf r}_{-N}^{\overline{\alpha}_{-N}}\dots {\bf r}_{M-1}^{\overline{\alpha}_{M-1}})|^2=p^{Ns}.
 \end{equation}
 and for any vector
 $({\bf a}_{-1},{\bf a}_{-2},\dots,{\bf a}_{-N}) \ne(0,0,\dots,0),\ ({\bf a}_j\in GF(p^s))$
  \small
 \begin{equation}                                                                     \label{eq6.3}
  \sum_{\overline{\alpha}_{-1},\dots,\overline{\alpha}_{-N}}\exp\left(\frac{2\pi
 i}{p}({\bf a}^{(0)}_{-1}\alpha_{-1}^{(0)}+\dots +{\bf
 a}^{(s-1)}_{-1}\alpha_{-1}^{(s-1)}+\dots+
 {\bf a}_{-N}^{(0)}\alpha_{-N}^{(0)}+\dots+{\bf a}_{-N}^{(s-1)}\alpha_{-N}^{(s-1)})\right)\times
 $$
  $$
  \times\sum_{\overline{\alpha}_0,\overline{\alpha}_1,\dots,\overline{\alpha}_{M-1}}|\hat\varphi((K^+_{-N})^\bot
  {\bf r}_{-N}^{\overline{\alpha}_{-N}}\dots {\bf r}_{M-1}^{\overline{\alpha}_{M-1}})|^2=0
 \end{equation}
 \large

 Let $(\varphi(x\dot-h))_{h\in H_0}$  be an orthonormal system. Using  the Plansherel equality and
 Lemma  5.1 we have
 $$
 \delta_{h_1h_2}=\int\limits_{K^+}\varphi(x\dot-h_1)\overline{\varphi(x\dot-h_2)}\,d\mu(x)=
 \int\limits_{ (K^+_M)^\bot}|\hat\varphi(\chi)|^2(\chi,h_2\dot-
 h_1)\,d\nu(\chi)=
 $$
 $$
 =\sum_{\overline{\alpha}_{-N},\dots,\overline{\alpha}_0,\dots,\overline{\alpha}_{M-1}}\int\limits_{(K^+_{-N})^\bot
 {\bf r}_{-N}^{\overline{\alpha}_{-N}}\dots {\bf r}_0^{\overline{\alpha}_0}\dots
 {\bf r}_{M-1}^{\overline{\alpha}_{M-1}}}|\hat\varphi(\chi)|^2(\chi,h_2\dot-
 h_1)\,d\nu(\chi)=
 $$
 $$
 =\sum_{\overline{\alpha}_{-N},\dots,\overline{\alpha}_{M-1}}|\hat\varphi((K^+_{-N})^\bot
 {\bf r}_{-N}^{\overline{\alpha}_{-N}}\dots
 {\bf r}_{M-1}^{\overline{\alpha}_{M-1}}|^2\int\limits_{(K^+_{-N})^\bot
 {\bf r}_{-N}^{\overline{\alpha}_{-N}}\dots
 {\bf r}_{M-1}^{\overline{\alpha}_{M-1}}}(\chi,h_2\dot- h_1)\,d\nu(\chi)=
 $$
 $$
 =p^{-sN}{\bf
 1}_{K^+_{-N}}(h_2\dot-h_1)\times
 $$
 $$
 \times\sum_{\overline{\alpha}_{-N},\dots,\overline{\alpha}_{M-1}}|\hat\varphi((K^+_{-N})^\bot
 {\bf r}_{-N}^{\overline{\alpha}_{-N}}\dots {\bf r}_0^{\overline{\alpha}_0}\dots
 {\bf r}_{M-1}^{\overline{\alpha}_{M-1}}|^2({\bf r}_{-N}^{\overline{\alpha}_{-N}}\dots
 {\bf r}_0^{\overline{\alpha}_0}\dots {\bf r}_{M-1}^{\overline{\alpha}_{M-1}},h_2\dot-h_1).
 $$
 If  $h_2=h_1$, we obtain the equality (\ref{eq6.2}). If $h_2\ne
 h_1$ then
 \begin{equation}                                                        \label{eq6.4}
 h_2\dot-h_1={\bf a}_{-1}g_{-1}\dot+\dots\dot+{\bf a}_{-N}g_{-N}\in K^+_{-N}
 \end{equation}
 or
 \begin{equation}                                                      \label{eq6.5}
 h_2\dot-h_1={\bf a}_{-1}g_{-1}\dot+\dots\dot+{\bf a}_{-N}g_{-N}\dot+
 \dots\dot+{\bf a}_{-l}g_{-l}\in K^+\backslash K^+_{-N}.
  \end{equation}
  If the condition (\ref{eq6.5}) are fulfilled,  then ${\bf 1}_{(K_{-N}^+)}(h_2\dot-h_1)=0$.
  If the condition  (\ref{eq6.4})
 are fulfilled, then
 $$
 {\bf 1}_{K^+_{-N}}(h_2\dot-h_1)=1,
 $$
 $$
  ({\bf r}_{-N}^{\overline{\alpha}_{-N}}\dots
 {\bf r}_0^{\overline{\alpha}_0}\dots
 {\bf r}_{M-1}^{\overline{\alpha}_{M-1}},h_2\dot-h_1)=
 ({\bf r}_{-N}^{{\bf a}_{-N}},\overline{\alpha}_{-N}g_{-N})
 \dots
 ({\bf r}_{-1}^{{\bf a}_{-1}},\overline{\alpha}_{-1}g_{-1}).
 $$
 Using the equality $({\bf r}_k^{{\bf a}_k},{\bf u}g_k)=\prod\limits_{l=0}^{s-1}e^{\frac{2\pi
i}{p}u^{(l)}a_k^{(l)}}$  we obtain the equality
 (\ref{eq6.3}).The conversely may be proved by analogy.

 Let as show now if for any vector  $({\bf a}_{-1},{\bf a}_{-2},\dots,{\bf a}_{-N})\ne(0,0,\dots,0)$
 the conditions   (\ref{eq6.2})  (\ref{eq6.3}) are fulfilled, then for any
   $\overline{\alpha}_{-N},\overline{\alpha}_{-N+1},\dots,\overline{\alpha}_{-1}\in FG(p^s)$
 \begin{equation}                                                                   \label{eq6.6}
  \sum_{\overline{\alpha}_{0},\overline{\alpha}_1,\dots,\overline{\alpha}_{M-1}}|\hat\varphi((K^+_{-N})^\bot
  {\bf r}_{-N}^{\overline{\alpha}_{-N}}\dots {\bf r}_0^{\overline{\alpha}_0}\dots
  {\bf r}_{M-1}^{\overline{\alpha}_{M-1}})|^2=1.
 \end{equation}
 Let us denote
 $$
  n=\sum_{j=1}^N\sum_{l=0}^{s-1} a_{-j}^{(l)}p^{s(j-1)},
  \;\;k=\sum_{j=1}^N \sum_{l=0}^{s-1}\alpha_{-j}^{(l)}p^{s(j-1)},
 $$
 $$
  \;\;C_{n,k}=\exp{\frac{2\pi
  i}{p}\left(\sum_{j=1}^N\sum_{l=0}^{s-1}\alpha_{-j}^{(l)}{a}_{-j}^{(l)}\right)}.
 $$
 and write the equalities  (\ref{eq6.2}) and (\ref{eq6.3}) as the system
 \begin{equation}   \label{eq6.7}
 \begin{array}{l}
  C_{0,0}x_0+C_{0,1}x_1+\dots+C_{0,p^{sN}-1}x_{p^{sN}-1}=p^{sN} \\
  C_{1,0}x_0+C_{1,1}x_1+\dots+C_{1,p^{sN}-1}x_{p^{sN}-1}=0 \\
  \dots \dots  \dots\dots\dots\dots\dots\dots\dots \dots \dots \dots \dots \dots\\
  C_{p^{sN}-1,0}x_0+C_{p^{sN}-1,1}x_1+\dots+C_{p^{sN}-1,p^{sN}-1}x_{p^{sN}-1}=0 \\
 \end{array}
 \end{equation}
  with unknowns
 $$
 x_k=\sum_{\overline{\alpha}_{0},\overline{\alpha}_1,\dots,\overline{\alpha}_{M-1}}|\hat\varphi((K^+_{-N})^\bot
 {\bf r}_{-N}^{\overline{\alpha}_{-N}}\dots {\bf r}_0^{\overline{\alpha}_0}\dots
 {\bf r}_{M-1}^{\overline{\alpha}_{M-1}})|^2.
 $$
 The matrix $(C_{n,k})$ is orthogonal. Indeed, if\\
 $({\bf a}_{-1},{\bf a}_{-2},\dots,{\bf a}_{-N})\ne({\bf a}'_{-1},{\bf a}'_{-2},\dots,{\bf a}'_{-N})$,
 i.e., $n\ne n'$ we obtain
 \normalsize
 $$
 \sum_{k=0}^{p^N-1}C_{n,k}\overline{C_{n',k}}=
 \sum_{\overline{\alpha}_{-1},\dots,\overline{\alpha}_{-N}}\exp\left(\frac{2\pi
 i}{p}\sum_{j=1}^N\sum_{l=0}^{s-1}  (a^{(l)}_{-j}-a'^{(l)}_{-j})\alpha_{-j}^{(l)}\right)=0,
 $$
 \large
 so at least one of differences  ${\bf a}_{-l}-{\bf a}'_{-l}\ne 0$.
 So, the system (\ref{eq6.7}) has unique solution.
 It is evident that  $x_k=1$ is  a solution of this system. This
 means that
 (\ref{eq6.6}) is fulfil, and the necessity is proved. The
 sufficiency is evident. $\square$

 Now we obtain a necessary and sufficient conditions for function  $m_0(\chi)$ to be
 a mask on the class
 $\mathfrak D_{-N}((K^+_M)^\bot)$, i.e. there exists
 $\hat\varphi\in\mathfrak D_{-N}((K^+_M)^\bot)$ for which
 \begin{equation}                                                                    \label{eq6.8}
 \hat\varphi(\chi)=m_0(\chi)\hat\varphi(\chi{\cal A}^{-1}).
 \end{equation}
 If $m_0(\chi)$ is a mask of (\ref{eq6.8}) then\\
 P1) $m_0(\chi)$
 is constant on cosets  $(K^+_{-N})^\bot\zeta$,\\
 P2) $m_0(\chi)$ is periodic with any period
  ${\bf r}_1^{\overline{\alpha}_1}{\bf r}_2^{\overline{\alpha}_2}\dots
 {\bf r}_l^{\overline{\alpha}_l}$, $\overline{\alpha}_j\in GF(p^s)$, \\
 P3)
 $m_0((K^+_{-N})^\bot)=1$. \\
 Therefore we will assume that  $m_0$
 satisfies these conditions.
 Let
 $$
 E_k\subset (K^+_k)^\bot\setminus (K^+_{k-1})^\bot\;\ ,(k=-N+1,-N+2,\dots,0,1,\dots,M,M+1)
 $$
 be a set,  on which $m_0(E_k)=0$.
 Since  $m_0(\chi)$ is constant on cosets
 $(K^+_{-N})^\bot\zeta$, it follows that $E_k$ is a union of such  cosets or $E_k=\emptyset$.

  \begin{theorem}
  $m_0(\chi)$ is a mask of some equation on the class $\mathfrak
  D_{-N}((K^+_M)^\bot)$ if and only if
  \begin{equation}                                                              \label{eq6.9}
   m_0(\chi)m_0(\chi{\cal A}^{-1})\dots m_0(\chi{\cal A}^{-M-N})=0
  \end{equation}
  on  $(K^+_{M+1})^\bot\setminus (K^+_M)^\bot$.
 \end{theorem}
 {\bf Proof.}  Indeed, if
 (\ref{eq6.9}) is true  we set
 $$
 \hat\varphi(\chi)=\prod\limits_{k=0}^\infty m_0(\chi{\cal
 A}^{-k})\in \mathfrak D_{-N}( (K^+_M)^\bot).
 $$
 Then $\hat\varphi(\chi)=m_0(\chi)\hat\varphi(\chi{\cal
 A}^{-1}) $  and
 $$
 m_0(\chi)=\sum_{h\in H_0^{(N+1)}}\beta_h\overline{(\chi{\cal
 A}^{-1},h)}
 $$
 for some  $\beta_h$. Therefore $m_0(\chi)$ is a mask. Inversely
  let $m_0(\chi)$ be a mask, i.e.
 $\hat\varphi(\chi)=m_0(\chi)\hat\varphi(\chi{\cal A}^{-1})\in \mathfrak D_{-N}((K^+_M)^\bot)$.
  From it we find
 $$
 \hat\varphi(\chi)=m_0(\chi)m_0(\chi{\cal A}^{-1})\dots
 m_0(\chi{\cal A}^{-M-N})\hat\varphi(\chi{\cal A}^{-M-N-1}),
 $$
 and  $\hat\varphi(\chi{\cal A}^{-M-N-1})=1$ on $(K^+_{M+1})^\bot$.
 Since  $\hat\varphi(\chi)=0$ on $(K^+_{M+1})^\bot\setminus (K^+_M)^\bot$,
 it follows
 $$
 m_0(\chi)m_0(\chi{\cal A}^{-1})\dots m_0(\chi{\cal A}^{-M-N})=0
 $$
 on  $(K^+_{M+1})^\bot\setminus (K^+_M)^\bot$. $\square$

\begin{lemma}
 Let $\hat\varphi\in \mathfrak D_{-N}((K^+_M)^\bot)$ be a solution
 of the refinement equation
 $$
 \hat\varphi(\chi)=m_0(\chi)\hat\varphi(\chi{\cal A}^{-1})
 $$
 and $(\varphi(x\dot-h))_{h\in H_0}$ be an orthonormal system.
 Then for any \\
 $\overline{\alpha}_{-N},\alpha_{-N+1},\dots,\overline{\alpha}_{-1}\in GF(p^s)$
 \begin{equation}                                                                           \label{eq6.10}
 \sum_{\overline{\alpha}_0\in GF(p^s)}|m_0((K^+_{-N})^\bot {\bf
 r}_{-N}^{\overline{\alpha}_{-N}}{\bf r}_{-N+1}^{\overline{\alpha}_{-N+1}}\dots {\bf
 r}_{-1}^{\overline{\alpha}_{-1}}{\bf r}_{0}^{\overline{\alpha}_{0}})|^2=1.
 \end{equation}
 \end{lemma}
 {\bf Proof.} Since $\hat\varphi\in\mathfrak D_{-N}((K^+_M)^\bot)$,
 it follows that $\hat\varphi((K^+_{M+1})^\bot\setminus
 ((K^+_M)^\bot)=0$. Using  theorem 6.1 we have
 $$
 1=\sum_{\overline{\alpha}_0,\overline{\alpha}_1,\dots, \overline{\alpha}_{M-1}\in
 GF(p^s)}|\hat\varphi((K^+_{-N})^\bot {\bf
 r}_{-N}^{\overline{\alpha}_{-N}}\dots {\bf r}_0^{\overline{\alpha}_0}\dots {\bf
 r}_{M-1}^{\overline{\alpha}_{M-1}})|^2=
 $$
 \small
 $$
 =\sum_{\overline{\alpha}_0,\dots, \overline{\alpha}_{M-1},\overline{\alpha}_{M}\in
 GF(p^s)}|\hat\varphi((K^+_{-N})^\bot {\bf
 r}_{-N}^{\overline{\alpha}_{-N}}\dots {\bf r}_0^{\overline{\alpha}_0} \dots {\bf
 r}_{M-1}^{\overline{\alpha}_{M-1}}{\bf
 r}_{M}^{\overline{\alpha}_{M}})|^2=\sum_{\overline{\alpha}_0=0}^{p-1}|m_0((K^+_{-N})^\bot
 {\bf r}_{-N}^{\overline{\alpha}_{-N}}\dots  {\bf r}_{0}^{\overline{\alpha}_{0}})|^2
 $$
 \large
 $$
 \cdot\sum_{\overline{\alpha}_1,\dots, \overline{\alpha}_{M-1},\overline{\alpha}_M\in
 GF(p^s)}|\hat\varphi((K^+_{-N})^\bot {\bf
 r}_{-N}^{\overline{\alpha}_{-N}+1}\dots {\bf r}_{-1}^{\overline{\alpha}_{0}} {\bf
 r}_0^{\overline{\alpha}_1}\dots {\bf r}_{M-2}^{\overline{\alpha}_{M-1}}
 {\bf r}_{M-1}^{\overline{\alpha}_{M}})|^2=
 $$
 $$
 =\sum_{\overline{\alpha}_0\in GF(p^s)}|m_0((K_{-N}^+)^\bot
 {\bf r}_{-N}^{\overline{\alpha}_{-N}}\dots  {\bf r}_0^{\overline{\alpha}_0})|^2.\;\;\square
 $$

\begin{theorem}
Suppose the function $m_0(\chi)$ satisfies the conditions
P1,P2,P3, (\ref{eq6.9}), and the function
$$
\hat\varphi(\chi)=\prod\limits_{n=0}^\infty m_0(\chi{\cal A}^{-n})
$$
satisfies the condition  (\ref{eq6.1}). Then $\varphi\in \mathfrak
D_M(K^+_{-N})$  and generates an orthogonal MRA.
\end{theorem}
 {\bf Proof.} It is evident that
  $\hat \varphi\in \mathfrak  D_{-N}((K^+_{M})^\bot)$,
 $\hat\varphi(\chi)=m_0(\chi)\hat\varphi(\chi{\cal A}^{-1})$, and
 $(\varphi(x\dot-h))_{h\in H_0}$ is an orthonormal system. From
 theorems 5.1, 5.2, 5.3 we find that the function $\varphi$
 generates an orthogonal MRA. $\square$

 \section{(N,M)-elementary sets}

 So,  to find a refinable function that generates orthogonal MRA,
 we need take  a function $m_0(\chi)$ that satisfies conditions P1,
 P2, P3, (\ref{eq6.9}), construct the function
 $$
 \hat\varphi(\chi)=\prod\limits_{k=0}^\infty m_0(\chi{\cal
 A}^{-k})\in \mathfrak D_{-N}((K_M^+)^\bot)
 $$
 and check that the system $\varphi(x\dot-h)_{h\in H_0}$ is
 orthonormal. We want to give a simple condition under witch
  the system of shifts $\varphi(x\dot-h)_{h\in H_0}$ is
 orthonormal.

 \begin{definition}
 Let $N,M\in \mathbb N$. A set $E \subset X$ is called
 $(N,M)$-elementary if $E$ is disjoint union of $p^{sN}$ cosets
 $$
  (K_{-N}^+)^\bot\zeta_j=
 (K_{-N}^+)^\bot\underbrace{{\bf r}_{-N}^{\overline{\alpha}_{-N}}{\bf r}_{-N+1}^{\overline{\alpha}_{-N+1}}\dots
 {\bf r}_{-1}^{\overline{\alpha}_{-1}}}_{\xi_j}\underbrace{{\bf r}_{0}^{\overline{\alpha}_{0}}\dots
 {\bf r}_{M-1}^{\overline{\alpha}_{M-1}}}_{\eta_j}=(K_{-N}^+)^\bot\xi_j\eta_j,
 $$
 $j=0,1,...,p^{sN}-1,
 j=\sum_{l=0}^{N-1}(\alpha_{-N+l}^{(0)}+\alpha_{-N+l}^{(1)}p+\dots+\alpha_{-N+l}^{(s-1)}p^{s-1})p^{sl}$
 $(\overline{\alpha}_{\nu}\in GF(p^s))$ such that\\
 1) $\bigsqcup\limits_{j=0}^{p^{sN}-1}
 (K_{-N}^+)^\bot\xi_j=(K_{0}^+)^\bot$, $(K_{-N}^+)^\bot\xi_0=(K_{-N}^+)^\bot$,\\
 2) for any $l=\overline{0,M+N-1}$ the intersection $(
 (K_{-N+l+1}^+)^\bot\setminus (K_{-N+l}^+)^\bot)\bigcap
 E\ne\emptyset$.
 \end{definition}
 \begin{lemma}
 The set $H_0\subset  K$ is an orthonormal  system on any
 $(N,M)$-elementary set $E\subset X$.
 \end{lemma}
 {\bf Proof.} Let $h,g\in H_0$. Using the definition of $(N,M)$-elementary set we
 have
 $$
 \int\limits_E(\chi,h)\overline{(\chi,g)}\,d\nu(\chi)=\sum_{j=0}^{p^{sN}-1}\int\limits_{
 (K_{-N}^+)^\bot\zeta_j}(\chi,h)\overline{(\chi,g)}\,d\nu(\chi)=
 $$
 $$
  =\sum\limits_{j=0}^{p^{sN}-1}\int\limits_X{\bf 1}_{
  (K_{-N}^+)^\bot\zeta_j}(\chi)(\chi,h)\overline{(\chi,g)}\,d\nu(\chi)=
 $$
 $$
  =\sum\limits_{j=0}^{p^{sN}-1}\int\limits_X{\bf 1}_{
  (K_{-N}^+)^\bot\zeta_j}(\chi\eta_j)(\chi\eta_j,h)\overline{(\chi\eta_j,g)}\,d\nu(\chi)=
 $$
 $$
 =\sum\limits_{j=0}^{p^{sN}-1}\int\limits_X{\bf 1}_{
 (K_{-N}^+)^\bot\xi_j}(\chi)(\chi,h)\overline{(\chi,g)}(\eta_j,h)\overline{(\eta_j,g)}\,d\nu(\chi).
 $$
 Since
 $$
 (\eta_j,h)=({\bf r}_0^{\overline{\alpha}_0}{\bf r}_1^{\overline{\alpha}_1}\dots
 {\bf r}_{M-1}^{\overline{\alpha}_{M-1}},{\bf a}_{-1}g_{-1}\dot+
 {\bf a}_{-2}g_{-2}\dot+\dots\dot+{\bf a}_{-l}g_{-l})=1,
 $$
 $$
 (\eta_j,g)=({\bf r}_0^{\overline{\alpha}_0}{\bf r}_1^{\overline{\alpha}_1}\dots
 {\bf r}_{M-1}^{\overline{\alpha}_{M-1}},{\bf b}_{-1}g_{-1}\dot+
 {\bf b}_{-2}g_{-2}\dot+\dots\dot+{\bf b}_{-l}g_{-l})=1,
 $$
 then
 $$
 \int\limits_E(\chi,h)\overline{(\chi,g)}\,d\nu(\chi)=
 \sum\limits_{j=0}^{p^{sN}-1}\int\limits_{(K_{-N}^+)^\bot\xi_j}(\chi,h)\overline{(\chi,g)}\,d\nu(\chi)=
 $$
 $$
 \int\limits_{(K_{0}^+)^\bot}(\chi,h)\overline{(\chi,g)}\,d\nu(\chi)
 =\delta_{h,g}$$. {$\square$}
 \begin{theorem}
 Let
  $K=F^{(s)}$ be a local field with positive characteristic
 $p$,  $E\subset (K_M^+)^\bot$ an $(N,M)$-elementary set. If\
 $|\hat\varphi(\chi)|={\bf 1}_E(\chi)$ on $X$ then the system of
 shifts $(\varphi(x\dot-h))_{h\in H_0}$ is an orthonormal system on
 $K$.
 \end{theorem}
 {\bf Proof.}
 Let $\tilde H_0\subset H_0$ be a finite set. Using
 the Plansherel equation we have
 $$
  \int\limits_{K}\varphi(x\dot-g)\overline{\varphi(x\dot-g)}\,d\mu(x)=
  \int\limits_X|\hat\varphi(\chi)|^2\overline{(\chi,g)}(\chi,h)d\nu(\chi)=
  \int\limits_E(\chi,h)\overline{(\chi,g)}d\nu(\chi)=
 $$
 $$
  =\sum_{j=0}^{p^{sN}-1}\int\limits_{(K_{-N}^+)^\bot\zeta_j}(\chi,h)\overline{(\chi,g)}\,d\nu(\chi).
 $$
 Transform the inner integral
 $$
  \int\limits_{(K_{-N}^+)^\bot\zeta_j}(\chi,h)\overline{(\chi,g)}\,d\nu(\chi)=\int\limits_X{\bf
  1}_{(K_{-N}^+)^\bot\zeta_j}(\chi)(\chi,h)\overline{(\chi,g)}\,d\nu(\chi)=
 $$
 $$
  =\int\limits_X{\bf 1}_{(K_{-N}^+)^\bot\zeta_j}(\chi\eta_j)(\chi\eta_j,h\dot-g)\,d\nu(\chi)=
  \int\limits_X{\bf 1}_{(K_{-N}^+)^\bot\xi_j}(\chi)(\chi\eta_j,h\dot-g)\,d\nu(\chi)=
 $$
 $$
  =\int\limits_{(K_{-N}^+)^\bot\xi_j}(\chi\eta_j,h\dot-g)\,d\nu(\chi).
 $$
 Repeating the arguments of lemma 7.1 we obtain
 $$
  \int\limits_{K}\varphi(x\dot-h)\overline{\varphi(x\dot-g)}\,d\mu(x)=\delta_{h,g}.\;\;\square
 $$

 \section{Trees and wavelets} Let $K=F^{(s)}$ be a local field of characteristic $p$.    In this section we reduce the problem of construction of step
refinable function on the field $K$ to construction of some tree.

 We will consider some
special class of refinable functions $\varphi(\chi)$ for which
$|\hat\varphi(\chi)|$ is a characteristic function of a set.
Define this class.
\begin{definition}
A mask $m_0(\chi)$ is called $N$-elementary $(N\in\mathbb N_0)$ if
$m_0(\chi)$ is constant on cosets $(K_{-N}^+)^\bot\chi$, its
modulus $m_0(\chi)$ has two values only: 0 and 1, and
$m_0((K_{-N}^+)^\bot)=1$. The refinable function $\varphi(x)$ with
Fourier transform
$$
\hat\varphi(\chi)=\prod\limits_{n=0}^\infty m_0(\chi{\cal A}^{-n})
$$
is called $N$-elementary too. $N$-elementary function $\varphi$
 is called $(N,M)$-elemen\-tary if  $\hat\varphi(\chi)\in \mathfrak D_{-N}(K_{M}^\bot)$.
  In this case  we will
 call the Fourier transform $\hat\varphi(\chi)$  $(N,M)$-elementary, also.
\end{definition}

 \begin{definition}
 Let $\tilde
 E=\bigsqcup\limits_{\overline{\alpha}_{-1},\overline{\alpha}_0}(K_{-1}^+)^\bot
 {\bf r}_{-1}^{\overline{\alpha}_{-1}}{\bf
 r}_0^{\overline{\alpha}_{0}}\subset (K_1^+)^\bot$ be an
 $(1,1)$-elementary set. We say that the set $\tilde E_X$ is a
 periodic extension of $\tilde E$ if
 $$
 \tilde E_X=\bigcup
 \limits_{l=1}^\infty\bigsqcup\limits_{\overline{\alpha}_1,\dots,\overline{\alpha}_l\in
 GF(p^s)}\tilde E {\bf r}_1^{\overline{\alpha}_1}{\bf
 r}_2^{\overline{\alpha}_2}\dots {\bf r}_l^{\overline{\alpha}_l}.
 $$
We say that the set $\tilde E$ generates an $(1,M)$ elementary set
$E$, if $\bigcap\limits_{n=0}^\infty \tilde E_X{\cal A}^n=E$.
\end{definition}


 Let us write the set $GF(p^s)$ in the form
 $$
 \{{\bf 0},{\bf u}_1,{\bf u}_2,\dots,{\bf u}_q,\overline{\alpha}_1,\overline{\alpha}_2,
 \dots,\overline{\alpha}_{p^s-q-1}\}=V,\;\;{\bf 0}={\bf u}_0,
 $$
 where $1\le q\le p^s-1$. We will consider the set $V$ as a set of
 vertices. By
 $T({\bf 0},{\bf u}_1,{\bf u}_2,\dots,{\bf u}_q,\overline{\alpha}_1,\overline{\alpha}_2,
 \dots,\overline{\alpha}_{p^s-q-1})=T(V)$
 we will denote a rooted tree on the set of vertices $V$, where ${\bf 0 }$ is a root,
 ${\bf u}_1,{\bf u}_2,\dots,{\bf u}_q$ are first level vertices,
 $\overline{\alpha}_1,\overline{\alpha}_2,\dots,\overline{\alpha}_{p^s-q-1}$ are remaining
 vertices.\\
 For example for $p=3, s=2, q=2, {\bf u}_1=(2,1), {\bf u}_2=(1,0)$
  we have trees

\begin{picture}(70,60)
    \put(24,8){$(0,0)$}

   \put(28,12){\vector(-1,1){6}}
   \put(14,20){$(2,1)$}

  \put(32,12){\vector(1,1){6}}
  \put(33,20){$(1,0)$}

  \put(16,26){\vector(-1,1){6}}
  \put(4,33){$(1,2)$}

   \put(22,26){\vector(1,1){6}}
   \put(28,33){$(0,2)$}

  \put(44,26){\vector(1,1){6}}
  \put(48,33){$(2,2)$}

  \put(34,40){\vector(1,1){6}}
  \put(34,48){$(1,1)$}

  \put(14,40){\vector(1,1){6}}
  \put(14,48){$(2,0)$}
\end{picture}
or
\begin{picture}(70,60)
    \put(24,8){$(0,0)$}

   \put(28,12){\vector(-1,1){6}}
   \put(14,20){$(2,1)$}

  \put(32,12){\vector(1,1){6}}
  \put(33,20){$(1,0)$}

  \put(16,26){\vector(-1,1){6}}
  \put(4,33){$(1,2)$}

   \put(22,26){\vector(1,1){6}}
   \put(24,33){$(2,2)$}

  \put(44,26){\vector(1,1){6}}
  \put(48,33){$(0,2)$}

  \put(34,40){\vector(1,1){6}}
  \put(34,48){$(1,1)$}

  \put(56,40){\vector(1,1){6}}
  \put(56,48){$(2,0)$}
\end{picture}

\hspace*{15mm} Figure 1 \hspace*{40mm} Figure 2\\
 and so on.

 For any tree path
 $$
 P_j=({\bf
 0}\rightarrow{\bf u}_j\rightarrow\overline{\alpha}_{l-1}\rightarrow\overline{\alpha}_{l-2}\rightarrow\dots \rightarrow
 \overline{\alpha}_{0}\rightarrow\overline{\alpha}_{-1})
 $$
 we construct the set of cosets
 \begin{equation}              \label{eq8.1}
   (K_{-1}^+)^\bot {\bf r}_{-1}^{{\bf u}_j},(K_{-1}^+)^\bot
    {\bf r}_{-1}^{\overline{\alpha}_{l-1}}{\bf r}_{0}^{{\bf u}_j},(K_{-1}^+)^\bot
    {\bf r}_{-1}^{\overline{\alpha}_{l-2}}{\bf r}_{0}^{\overline{\alpha}_{l-1}},\dots,
    (K_{-1}^+)^\bot {\bf r}_{-1}^{\overline{\alpha}_{0}}{\bf r}_{0}^{\overline{\alpha}_1},
    (K_{-1}^+)^\bot {\bf r}_{-1}^{\overline{\alpha}_{-1}}{\bf r}_{0}^{\overline{\alpha}_0}.
 \end{equation}
 For example for the tree from Figure 2 and the path
 $$
 (0,0)\rightarrow (1,0)\rightarrow (0,2)\rightarrow (2,0)
 $$
 we have 3 cosets
 $$
   (K_{-1}^+)^\bot {\bf r}_{-1}^{(1,0)},(K_{-1}^+)^\bot {\bf r}_{-1}^{(0,2)}
    {\bf r}_0^{(1,0)}, (K_{-1}^+)^\bot {\bf r}_{-1}^{(2,0)}{\bf r}_0^{(0,2)},
 $$
 for the path $(0,0), (2,1), (2,2)$ we have two cosets
 $$
   (K_{-1}^+)^\bot {\bf r}_{-1}^{(2,1)}, (K_{-1}^+)^\bot {\bf r}_{-1}^{(2,2)}{\bf r}_0^{(2,1)}.
 $$
We will represent the tree $T(V)$ as the tree
\begin{picture}(40,12)
  \put(14,-2){$(0,0)$}
   \put(18,2){\vector(-1,1){6}}
   \put(22,2){\vector(1,1){6}}
  \put(10,8){$T_1$}
  \put(17,9){\circle*{1}}
  \put(21,9){\circle*{1}}
  \put(25,9){\circle*{1}}
  \put(28,8){$T_q$}
  \end{picture}
where $T_j$ are tree branches of $T(V)$ with ${\bf u}_j$ as a
root.
 By $E_j$ denote a union of all cosets \eqref{eq8.1} for fixed $j$
 and set
 \begin{equation}\label{eq8.2}
\tilde E=\left(\bigsqcup\limits_{j=1}^q E_j\right)\bigsqcup
(K_{-1}^+)^\bot.
 \end{equation}
 It is clear that $\tilde E$ is an $(1,1)$ elementary set and $\tilde E\subset (K_{1}^+)^\bot$.
\begin{definition}
Let $\tilde E_X$ be a periodic extension of $\tilde E$. We say
that the tree $T(V)$ generates a set $E$, if
$E=\bigcap\limits_{n=0}^\infty\tilde E_X{\cal A}^n.$
\end{definition}

\begin{lemma}
Let $T(V)$ be a rooted tree with ${\bf 0}=(0,0,...,0)$ as a root.
Let $E\subset X$ be a set generated by the tree $T(V)$, $H$ a
hight of $T(V)$. Then $E$ is an $(1,H-2)$-elementary set.
\end{lemma}
{\bf Proof.} Let us denote
$$
m(\chi)={\bf 1}_{\tilde
E_X}(\chi),\;\;M(\chi)=\prod\limits_{n=0}^\infty m(\chi{\cal
A}^{-n}).
$$
First we note that $M(\chi)={\bf 1}_E(\chi)$. Indeed
$$
{\bf 1}_E(\chi)=1\Leftrightarrow\chi\in E\Leftrightarrow
\forall\,n,\;\chi{\cal A}^{-n}\in\tilde
E_X\Leftrightarrow\forall\,n,\;{\bf 1}_{\tilde E_X}(\chi{\cal
A}^{-n})=1\Leftrightarrow
$$
$$
\forall\,n,\;m(\chi{\cal
A}^{-n})=1\Leftrightarrow\prod\limits_{n=0}^\infty m(\chi{\cal
A}^{-n})=1\Leftrightarrow M(\chi)=1.
$$
 It means that $M(\chi)={\bf 1}_E(\chi)$.
 Now we will prove, that ${\bf 1}_E(\chi)=0$ for
 $\chi\in (K_{H-1}^+)^\bot\setminus (K_{H-2}^+)^\bot$.
 Since
 $\tilde  E_X\supset (K_{-1}^+)^\bot$
 it follows that
 ${\bf 1}_{\tilde E_X}((K_{H-1}^+)^\bot{\cal A}^{-H})={\bf
 1}_{\tilde E_X}((K_{-1}^+)^\bot)=1$.
 Consequently
 $$
  \prod\limits_{n=0}^\infty{\bf 1}_{\tilde E_X}(\chi{\cal
  A}^{-n})=\prod\limits_{n=0}^{H-1}{\bf 1}_{\tilde E_X}(\chi{\cal
  A}^{-n})
 $$
 for $\chi\in (K_{H-1}^+)^\bot\setminus (K_{H-2}^+)^\bot$. Let us
 denote $m((K_{-1}^+)^\bot {\bf r}_{-1}^{\bf i} {\bf r}_0^{\bf
 k})=\lambda_{{\bf i},{\bf k}}$. By the definition of cosets
 \eqref{eq8.1} $m((K_{-1}^+)^\bot {\bf r}_{-1}^{\bf i}{\bf
 r}_0^{\bf k})\ne 0\Leftrightarrow$ the pair $({\bf k},{\bf i})$ is
 an edge of the tree $T(V)$.

 We need prove that
 $$
  {\bf 1}_E((K_{-1}^+)^\bot {\bf
  r}_{-1}^{\overline{\alpha}_{-1}}{\bf
  r}_{0}^{\overline{\alpha}_{0}}\dots {\bf
  r}_{H-2}^{\overline{\alpha}_{H-2}})=0
 $$
 for $\overline{\alpha}_{H-2}\ne 0$. Since $\tilde E_X$ is a
 periodic extension of $\tilde E$ it follows that the function
 $m(\chi)={\bf 1}_{\tilde E_X}(\chi)$ is periodic with any period
 ${\bf r}_{1}^{\overline{\alpha}_{1}}{\bf
 r}_{2}^{\overline{\alpha}_{2}}\dots {\bf
 r}_{l}^{\overline{\alpha}_{l}}$, $l\in\mathbb N$, i.e. $m(\chi
 {\bf r}_{1}^{\overline{\alpha}_{1}}{\bf
 r}_{2}^{\overline{\alpha}_{2}}\dots {r\bf
 }_{l}^{\overline{\alpha}_{l}})=m(\chi)$ when $\chi\in
 (K_1^+)^\bot$. Using this fact we can write $M(\chi)$ for $\chi\in
  (K^+_{H-1})^\bot\setminus (K^+_{H-2})^\bot$ in the form
 $$
  M((K^+_{-1})^\bot\zeta)=M((K_{-1}^+)^\bot
  {\bf r}_{-1}^{\overline{\alpha}_{-1}}{\bf r}_{0}^{\overline{\alpha}_{0}}\dots
  {\bf r}_{H-2}^{\overline{\alpha}_{H-2}})=
$$
 $$
  =m((K_{-1}^+)^\bot {\bf r}_{-1}^{\overline{\alpha}_{-1}}{\bf
  r}_{0}^{\overline{\alpha}_{0}})m((K_{-1}^+)^\bot
  {\bf r}_{-1}^{\overline{\alpha}_{0}}{\bf r}_{0}^{\overline{\alpha}_{1}})
  \dots
 $$
 $$
  m((K_{-1}^+)^\bot {\bf r}_{-1}^{\overline{\alpha}_{H-3}}{\bf r}_{0}^{\overline{\alpha}_{H-2}})
  m((K_{-1}^+)^\bot {\bf r}_{-1}^{\overline{\alpha}_{H-2}})=
 $$
$$
=\lambda_{\overline{\alpha}_{-1},\overline{\alpha}_0}\lambda_{\overline{\alpha}_{0},\overline{\alpha}_1}
\dots\lambda_{\overline{\alpha}_{H-3},\overline{\alpha}_{H-2}}\lambda_{\overline{\alpha}_{H-2},{\bf
0}},\; \overline{\alpha}_{H-2}\ne {\bf 0}.
$$
 If $\lambda_{\overline{\alpha}_{H-2},{\bf 0}}=0$ then
 $M((K_{-1}^+)^\bot\zeta)=0$. Let $\lambda_{\overline{\alpha}_{H-2},{\bf 0}}\ne 0$. It means
 that $\overline{\alpha}_{H-2}={\bf u}_j$ for some $j=\overline{1,q}$. If
 $\lambda_{\overline{\alpha}_{H-3},\overline{\alpha}_{H-2}}=0$ then
 $M((K_{-1}^+)^\bot\zeta)=0$. Assume that
 $\lambda_{\overline{\alpha}_{H-3},\overline{\alpha}_{H-2}}\ne 0$. It is true iff the
 pair $(\overline{\alpha}_{H-2},\overline{\alpha}_{H-3})$ is an edge of $T(V)$.
  Repeating  these arguments, we obtain a path
 $$
 ({\bf 0}\rightarrow{\bf u}_j=\overline{\alpha}_{H-2}\rightarrow\overline{\alpha}_{H-3}
 \rightarrow\dots \rightarrow\overline{\alpha}_l)
 $$ of the tree
 $T(V)$. Since $ {\rm hight}(T)= H$ it follows that $l\ge 0$.
 Consequently $(\overline{\alpha}_l,\overline{\alpha}_{l-1})$ is not edge and
 $\lambda_{\overline{\alpha}_{l-1},\overline{\alpha}_l}=0$, where $l\ge 0$. It means
 that $M((K_{-1}^+)^\bot\zeta) =0$.

 Now we prove that $E$ is $(1,H-2)$-elementary set. Indeed, any
 path
 $$
 ({\bf 0}\rightarrow{\bf
 u}_j=\overline{\alpha}_{l-1}\rightarrow\overline{\alpha}_{l-2}\rightarrow\dots
 \rightarrow\overline{\alpha}_0\rightarrow\overline{\alpha}_{-1})
 $$
 defines the coset $(K_{-1}^+)^\bot
 {\bf r}_{-1}^{\overline{\alpha}_{-1}}{\bf r}_{0}^{\overline{\alpha}_{0}}\dots
 {\bf r}_{l-1}^{\overline{\alpha}_{l-1}}\subset E$. But for any
 $\overline{\alpha}_{-1}\in GF(p^s)$ there exists unique path with
 endpoint $\overline{\alpha}_{-1}$ and starting point zero. It means that $E$
 is $(1,H-2)$-elementary set. $\square$

 \begin{theorem}
 Let $M,s\in\mathbb N$, $p^s\ge 3$. Let $E\subset (K_M^+)^\bot$ be an
 $(1,M)$-elementary set, $\hat\varphi\in\mathfrak D_{-1}((K_M^+)^\bot)$,
 $|\hat\varphi(\chi)|={\bf 1}_E(\chi)$,
 $\hat\varphi(\chi)$ the solution of the equation
 \begin{equation}              \label{eq8.3}
 \hat\varphi(\chi)=m_0(\chi)\hat\varphi(\chi{\cal A}^{-1}),
 \end{equation}
 where $m_0(\chi)$ is an 1-elementary mask. Then there exists a
 rooted tree $T(V)$ with ${\rm height}(T)=M+2$ that generates the
 set $E$.
 \end{theorem}
 {\bf Prof.} Since the set $E$ is $(1,M)$-elementary set and
 $|\hat\varphi(\chi)|={\bf 1}_E(\chi)$, it follows from theorem 7.1
 that the system $(\varphi(x\dot-h))_{h\in H_0}$ is an orthonormal
 system in $L_2(K)$. Using the theorem 6.1 we obtain that
 for  any $\overline{\alpha}_{-1}\in GF(p^s)$
 \begin{equation}              \label{eq8.4}
  \sum_{\overline{\alpha}_0,\overline{\alpha}_1,\dots, \overline{\alpha}_{M-1}\in GF(p^s)}
  |\hat\varphi((K_{-1}^+)^\bot
  {\bf r}_{-1}^{\overline{\alpha}_{-1}}{\bf r}_{0}^{\overline{\alpha}_0}\dots
  {\bf r}_{M-1}^{\overline{\alpha}_{M-1}})|^2=1.
 \end{equation}
 Since $\hat\varphi$ is a solution of refinement equation
 \eqref{eq8.3} it follows from lemma 6.1 that for
 $\overline{\alpha}_{-1}\in GF(p^s)$
 \begin{equation}              \label{eq8.5}
    \sum_{\overline{\alpha}_0\in GF(p^s)}|m_0((K_{-1}^+)^\bot
    {\bf r}_{-1}^{\overline{\alpha}_{-1}}{\bf r}_{0}^{\overline{\alpha}_{0}})|^2=1.
 \end{equation}
 Let as denote $\lambda_{\overline{\alpha}_{-1},\overline{\alpha}_0}:=
 m_0((K_{-1}^+)^\bot {\bf r}_{-1}^{\overline{\alpha}_{-1}}{\bf r}_{0}^{\overline{\alpha}_{0}})$. Then we
 write \eqref{eq8.5} in the form
 \begin{equation}              \label{eq8.6}
 \sum_{\overline{\alpha}_0\in GF(p^s)}|\lambda_{\overline{\alpha}_{-1},\overline{\alpha}_0}|^2=1.
 \end{equation}
 Since the mask $m_0(\chi)$ is 1-elementary it follows that
 $|\lambda_{\overline{\alpha},\overline{\beta}}|$ take two values only: 0 or 1.

 Now we will construct the tree $T$. Let $\mathfrak U$ be a family
 of cosets $(K_{-1}^+)^\bot\zeta\subset (K_M^+)^\bot$ such that
 $\hat\varphi((K_{-1}^+)^\bot \zeta)\ne 0$ and
 $(K_{-1}^+)^\bot\notin \mathfrak U$. We can write a coset
 $(K_{-1}^+)^\bot\zeta\in \mathfrak U$ in the form
 $$
 (K_{-1}^+)^\bot {\bf r}_{-1}^{\overline{\alpha}_{-1}}{\bf
 r}_{0}^{\overline{\alpha}_{0}}\dots {\bf r}_{M-1}^{\overline{\alpha}_{M-1}}.
 $$
 If $(K_{-1}^+)^\bot\zeta\subset (K_{l}^+)^\bot\setminus
 (K_{l-1}^+)^\bot$ $(l\le M)$ then
  $$
 (K_{-1}^+)^\bot\zeta=(K_{-1}^+)^\bot
 {\bf r}_{-1}^{\overline{\alpha}_{-1}}{\bf r}_{0}^{\overline{\alpha}_{0}}\dots
 {\bf r}_{l-1}^{\overline{\alpha}_{l-1}},\ \overline{\alpha}_{l-1} \neq {\bf 0}.
 $$

 Let ${\bf u}\neq {\bf 0}$. By $T_{\bf u}$ we denote the set of vectors
 $({\bf u},\overline{\alpha}_{n-1},\dots,\overline{\alpha}_{0},\overline{\alpha}_{-1})$ for which
 $(K_{-1}^+)^\bot {\bf r}_{-1}^{\overline{\alpha}_{-1}}{\bf r}_{0}^{\overline{\alpha}_{0}}\dots
 {\bf r}_{n-1}^{\overline{\alpha}_{n-1}}{\bf r}_n^{\bf u}\in \mathfrak U$. We will name the
 vector\\ $({\bf u},\overline{\alpha}_{n-1},\dots,{\overline{\alpha}_{0}},\overline{\alpha}_{-1})$
 as a path
 too. So $T_{\bf u}$ is the set of pathes with starting point ${\bf u}$, for
 which $\hat\varphi((K_{-1}^+)^\bot
 {\bf r}_{-1}^{\overline{\alpha}_{-1}}{\bf r}_{0}^{\overline{\alpha}_{0}}\dots
 {\bf r}_{n-1}^{\overline{\alpha}_{n-1}}{\bf r}_n^{\bf u})\ne 0$.
 Denote (it follow from \eqref{eq8.5}), if
 $\hat\varphi((K_{-1}^+)^\bot
 {\bf r}_{-1}^{\overline{\alpha}_{-1}}{\bf r}_{0}^{\overline{\alpha}_{0}}\dots
 {\bf r}_{n-1}^{\overline{\alpha}_{n-1}}{\bf r}_n^{\bf u})\ne 0$
 then \\
 $\hat\varphi((K_{-1}^+)^\bot
 {\bf r}_{-1}^{\overline{\alpha}_{-1}}{\bf r}_{0}^{\overline{\alpha}_{0}}\dots
 {\bf r}_{n-1}^{\overline{\alpha}_{n-1}}{\bf r}_n^{\bf u}
 {\bf r}_{n+1}^{\overline{\alpha}_{n+1}}) = 0$ for any $\overline{\alpha}_{n+1}\neq {\bf
 0}$.
  We will show that $T_{\bf u}$ is a
 rooted tree with ${\bf u}$ as a root.

 1) All vertices $\overline{\alpha}_j,{\bf u}$ of the path
 $({\bf u},\overline{\alpha}_{n-1},\dots,\overline{\alpha}_{0},\overline{\alpha}_{-1})$ are pairwise
 distinct. Indeed
 $$
   \hat\varphi((K_{-1}^+)^\bot
   {\bf r}_{-1}^{\overline{\alpha}_{-1}}{\bf r}_{0}^{\overline{\alpha}_{0}}\dots
   {\bf r}_{n-1}^{\overline{\alpha}_{n-1}}{\bf r}_n^{\bf u})=\lambda_{\overline{\alpha}_{-1},
   \overline{\alpha}_0}
   \lambda_{\overline{\alpha}_{0},\overline{\alpha}_1}\dots
   \lambda_{\overline{\alpha}_{n-1},{\bf u}}\overline{\lambda}_{\bf u,0}\ne
   0,\;{\bf u}\ne0.
 $$
 If $\overline{\alpha}_{n-1}={\bf u}$ then $|\lambda_{\bf u,u}|=|\lambda_{\bf u,0}|=1$
 that
 contradicts the equation \eqref{eq8.5}.\\
 If $\overline{\alpha}_{n-1}=0$ then $|\lambda_{\bf 0,u}|=|\lambda_{\bf 0,0}|=1$
 that contradicts the equation \eqref{eq8.5} too. Consequently
 $\overline{\alpha}_{n-1}\notin  \{\bf 0,u\}$. By analogy we obtain
 that\\
 $\overline{\alpha}_i\notin
 \{ {\bf 0,u},\overline{\alpha}_{n-1},\dots,\overline{\alpha}_{i+2},\overline{\alpha}_{i+1}\}$.

 2) If two pathes $({\bf
 u},\overline{\alpha}_{k-1},\dots,\overline{\alpha}_{0},\overline{\alpha}_{-1})$
 and $({\bf u},\overline{\beta}_{l-1},\dots,\overline{\beta}_{0},\overline{\beta}_{-1})$ have the
 common subpath
 $({\bf u},\overline{\alpha}_{k-1},\dots,\overline{\alpha}_{k-j+1},\overline{\alpha}_{k-j})=({\bf
 u},\overline{\beta}_{l-1},\dots,\overline{\beta}_{l-j+1},\overline{\beta}_{l-j})$ and
 $\overline{\alpha}_{k-j-1}\ne \overline{\beta}_{l-j-1}$ then
 $\{\overline{\alpha}_{-1},\overline{\alpha}_{0},\dots,
 \overline{\alpha}_{k-j-1}\}\bigcap\{\overline{\beta}_{-1},\overline{\beta}_{0},\dots,
 \overline{\beta}_{l-j-1}\}=\emptyset$.
 Indeed, let
 $$
    \{\overline{\alpha}_{-1},\overline{\alpha}_{0},\dots,
    \overline{\alpha}_{k-j-1}\}\bigcap\{\overline{\beta}_{-1},\overline{\beta}_{0},
    \dots,\overline{\beta}_{l-j-1}\}\ne\emptyset.
 $$
  Then there exists
  ${\bf v}\in\{\overline{\alpha}_{-1},\overline{\alpha}_{0},\dots,
  \overline{\alpha}_{k-j-1}\}\bigcap\{\overline{\beta}_{-1},\overline{\beta}_{0},
  \dots,\overline{\beta}_{l-j-1}\}$.\\
  Assume that ${\bf v}\ne \overline{\alpha}_{k-j-1}$. Then ${\bf v}=\overline{\alpha}_\nu$, $-1\le\nu
  \le k-j-2$ and ${\bf v}=\overline{\beta}_\mu$, $-1\le\mu\le l-j-1$. It follows that
 $$
 ({\bf u}=\overline{\alpha}_k,\dots,\overline{\alpha}_{k-j},\overline{\alpha}_{k-j-1},
 \dots,\overline{\alpha}_{\nu+1},\overline{\alpha}_{\nu}=
 \overline{\beta}_{\mu},\overline{\beta}_{\mu-1},\dots,\overline{\beta}_{0},\overline{\beta}_{-1})\in
 T_{\bf u}
 $$
 $$
 ({\bf u}=\overline{\beta}_{l},\dots,\overline{\beta}_{l-j}=
 \overline{\alpha}_{k-j},\overline{\beta}_{l-j-1},\dots,
 \overline{\beta}_{\mu+1},\overline{\beta}_\mu,\overline{\beta}_{\mu-1},\dots,
 \overline{\beta}_0,\overline{\beta}_{-1})\in T_u.
 $$
 So we have two different pathes with the same sheet $\overline{\beta}_{-1}$.
 But this contradicts theorem 6.1. This means that $T_{\bf u}$ has no
 cycles, consequently $T_{\bf u}$ is a graph with ${\bf u}$ as a root.

 3) By analogy we can proof that different trees $T_{\bf u}$ an $T_{\bf v}$ has
 no common vertices. It follows that the graph
 $T=({\bf 0},T_{{\bf u}_1},\dots,T_{{\bf u}_q})$ is a tree with ${\bf 0}$ as a root.

 4) It is evident that this tree generates refinable function
 $\hat\varphi$ with a mask $m_0$. Show that ${\rm height}(T)=M+2$.
 Indeed, since $\hat\varphi\in  \mathfrak D_{-1}((K_{M}^+)^\bot)$ it follows that there exists a
 coset $(K_{-1}^+)^\bot r_{-1}^{\overline{\alpha}_{-1}}r_{0}^{\overline{\alpha}_{0}}
 \dots {\bf r}_{M-1}^{\overline{\alpha}_{M-1}}$, $\overline{\alpha}_{M-1}\ne 0$ for which
 $$
 |\hat\varphi((K_{-1}^+)^\bot {\bf r}_{-1}^{\overline{\alpha}_{-1}}{\bf r}_{0}^{\overline{\alpha}_{0}}
 \dots {\bf r}_{M-1}^{\overline{\alpha}_{M-1}})|=1.
 $$
 This coset generates a path $({\bf 0},\overline{\alpha}_{M-1}={\bf u},\overline{\alpha}_{M-2},\dots,
 \overline{\alpha}_0,\overline{\alpha}_{-1})$ of $T$. This path contain $M+2$ vertex. It
 means that ${\rm height}(T)\ge M+2$. On the other hand there isn't
 coset $(K_{-1}^+)^\bot\zeta\in \mathfrak U$ with condition
 $(K_{-1}^+)^\bot\zeta\subset (K_{M+1}^+)^\bot\setminus (K_{M}^+)^\bot$,
 consequently  there isn't path with $L>M+2$. So
 ${\rm height}(T)=M+2$. Since ${\rm supp}\,\hat\varphi(\chi)$ is
 $(1,M)$-elementary set, it follows that the set of all vertices of
 the tree $T$ is the set $GF(p^s)$. The theorem is
 proved.  $\square$

\begin{definition}
 Let $T(V)$ be a rooted tree with ${\bf 0}$ as a root,  $H$ a hight
 of $T(V)$, $V=GF(p^s)$. Using cosets \eqref{eq8.1} we define the
 mask $m_0(\chi)$ in the subgroup $(K_{1}^+)^\bot$ as follows:
 $m_0((K_{-1}^+)^\bot)=1, m_0((K_{-1}^+)^\bot
 {\bf r}_{-1}^{\bf i}{r}_0^{\bf j})=\lambda_{\bf i,j}$, $|\lambda_{\bf i,j}|=1$ when
 $(K_{-1}^+)^\bot {\bf r}_{-1}^{\bf i}{\bf r}_0^{\bf j}\subset \tilde E$,
 (q.v. \eqref{eq8.2}), $|\lambda_{\bf i,j}|=0$ when $(K_{-1}^+)^\bot
 {\bf r}_{-1}^{\bf i}{\bf r}_0^{\bf j}\subset (K_{1}^+)^\bot\setminus\tilde E$.
 Let us extend the mask $m_0(\chi)$ on the $X\setminus (K_{1}^+)^\bot$ periodically, i.e. $m_0(\chi
 {\bf r}_1^{\overline{\alpha}_1}{\bf r}_2^{\overline{\alpha}_2}\dots {\bf r}_l^{\overline{\alpha}_l})
 =m_0(\chi)$. Then
 we say that the tree $T(V)$ generates the mask $m_0(\chi)$. Set
 $\hat \varphi(\chi)=\prod\limits_{n=0}^\infty m_0(\chi{\cal
 A}^{-n})$.
 It follows from Lemma 8.1 that\\
 1) ${\rm
 supp}\,\hat\varphi(\chi)\subset (K_{H-2}^+)^\bot$,\\
 2)
 $\hat\varphi(\chi)$ is $(1,H-2)$-elementary function,\\
 3)
 $(\varphi(x\dot-h))_{h\in H_0}$ is an orthonormal system.\\ In
 this case we say that the tree $T(V)$ generates the refinable
 function $\varphi(x)$.
 \end{definition}

 \begin{theorem}
 Let $p\ge 2$ be a prime number, $ s\in \mathbb N, p^s\ge 3,$
 $$
  V=\{{\bf 0},{\bf u}_1,{\bf u}_2,\dots,{\bf u}_q,{\bf
  a}_1,{\bf a}_2,\dots,{\bf a}_{p^s-q-1}\}
 $$
 a set of vertices, $T(V)$ a rooted tree, $\bf 0$ the root,
 ${\bf u}_1,{\bf u}_2,\dots,{\bf u}_q$ a first level vertices. Let $H$ be are height
 of $T(V)$. By $\varphi(x)$ denote the function generated by the
 $T(V)$. Then $\varphi(x)$ generate an orthogonal MRA on $F^{(s)}$.
 \end{theorem}
 {\bf Proof.} Since  $T(V)$ generates the function  $\varphi$, it
 follows  that 1)$\hat\varphi\in\mathfrak
 D_{-1}(K_1^+)^\bot)$, 2)$\hat\varphi(\chi)$ is $(1,H-2)$
 elementary function, 3)$\hat\varphi(\chi)$ is a solution of
 refinable equation \eqref{eq8.3}, 4)$(\varphi(x\dot-h))_{h\in
 H_0}$ is an orthonormal system. From the theorem 6.3 it follows
 that $\varphi(x)$ generates an orthogonal MRA. $\square$

{\bf Remark.} Now we can give a simple  algorithm for constructing
 non-Haar refinable function $\varphi(x)$. Let $T(V)$ be a tree on
 the set $V=GF(p^s)$. Construct a finite sequence $(\lambda_{\bf
 i,j})_{{\bf i,j}\in GF(p^s)}$ as follows: $\lambda_{\bf 0,0}=1$,
 $|\lambda_{\bf i,j}|=1$ if the pair $(\bf j,i)$ is an edge of $
 T(f)$. For any vertex $\overline{\alpha}_{-1}$ we take the path
 $({\bf 0}=\overline{\alpha}_{l+1},{\bf u}_j=\overline{\alpha}_l,
 \overline{\alpha}_{l-1},\dots,\overline{\alpha}_0,\overline{\alpha}_{-1})$ and suppose
 $$
  \hat\varphi((K_{-1}^+)^\bot
  {\bf r}_{-1}^{\overline{\alpha}_{-1}}{\bf r}_{0}^{\overline{\alpha}_{0}}\dots
  {\bf r}_{l-1}^{\overline{\alpha}_{l-1}}{\bf r}_{l}^{\overline{\alpha}_{l}}{\bf r}_{l+1}^{\bf 0})=
  \lambda_{\overline{\alpha}_{-1},\overline{\alpha}_0}\cdot
  \lambda_{\overline{\alpha}_{0},\overline{\alpha}_1}\cdot\dots
  \cdot\lambda_{\overline{\alpha}_{l-1},\overline{\alpha}_l}\cdot\lambda_{\overline{\alpha}_l,{\bf 0}}.
 $$
 Otherwise we suppose $\hat\varphi((K_{-1}^+)^\bot\zeta)=0$.
 Then $\varphi$ generates an orthogonal MRA on the field
 $GF(p^s)$.

{\bf Example.}  Let $p=s=2$. For these values we have trees

\noindent
\begin{picture}(60,40)
    \put(24,8){$(0,0)$}

   \put(28,12){\vector(-2,1){12}}
   \put(3,20){$(0,1)$}

  \put(30,12){\vector(0,1){6}}
  \put(24,20){$(1,1)$}

  \put(32,12){\vector(2,1){12}}
  \put(43,20){$(1,0)$}
\end{picture}
\begin{picture}(70,40)
    \put(24,8){$(0,0)$}

   \put(28,12){\vector(-2,1){12}}
   \put(3,20){$(0,1)$}

  \put(54,26){\vector(2,1){12}}
  \put(64,34){$(1,1)$}

  \put(32,12){\vector(2,1){12}}
  \put(43,20){$(1,0)$}
\end{picture}
\begin{picture}(60,40)
    \put(30,8){$(0,0)$}

   \put(28,24){\vector(-2,1){12}}
   \put(9,32){$(0,1)$}

  \put(38,12){\vector(0,1){6}}
  \put(30,20){$(1,1)$}

  \put(44,24){\vector(2,1){12}}
  \put(49,32){$(1,0)$}
 \end{picture}

\hspace*{10mm} Figure 3 \hspace*{30mm} Figure 4 \hspace*{40mm} Figure 5\\
 and so on. For the tree on figure 5 we obtain $\hat\varphi (\chi)$ in
 the form\\
 $\hat\varphi(K^+_{-1})=1,$  $\hat\varphi((K^+_{-1})^\bot{\bf
 r}_{-1}^{(1,1)})=\lambda_{1,1}$, $\hat\varphi((K^+_{-1})^\bot{\bf
 r}_{-1}^{(0,1)}{\bf
 r}_{0}^{(1,1)})=\lambda_{0,1}$,\\ $\hat\varphi((K^+_{-1})^\bot{\bf
 r}_{-1}^{(1,0)}{\bf
 r}_{0}^{(1,1)})=\lambda_{1,0}$. $|\lambda_{i,j}|=1$ and $\hat\varphi((K^+_{-1})^\bot
 \zeta)=0$ otherwise.
Suppose for simplicity $\lambda_{i,j}=1$. Then we can
calculate the scaling function
$$
\varphi(x)=\int\limits_X\hat\varphi(\chi)(\chi,x)d\,\nu(\chi)=
\int\limits_{(K_{-1}^+)^\bot}(\chi,x)d\,\nu(\chi)+
\int\limits_{(K_{-1}^+)^\bot{\bf
r}_{-1}^{(1,1)}}(\chi,x)d\,\nu(\chi)+
$$
$$
+\int\limits_{(K_{-1}^+)^\bot{\bf r}_{-1}^{(0,1)}{\bf
r}_{0}^{(1,1)}}(\chi,x)d\,\nu(\chi)+\int\limits_{(K_{-1}^+)^\bot{\bf
r}_{-1}^{(1,0)}{\bf
r}_{0}^{(1,1)}}(\chi,x)d\,\nu(\chi)=2^{-2}({\bf 1}_{K_{-1}^+}(x)+
$$
$$
+{\bf r}_{-1}^{(1,1)}(x){\bf 1}_{K_{-1}^+}(x)+{\bf
r}_{-1}^{(0,1)}(x){\bf r}_{0}^{(1,1)}(x){\bf 1}_{K_{-1}^+}(x)+{\bf
r}_{-1}^{(1,0)}(x){\bf r}_{0}^{(1,1)}(x){\bf 1}_{K_{-1}^+}(x))=
{\bf 1}_E(x)
$$
 where
 $$
 E=K^+_1 \bigsqcup (K^+_1\dot+(0,0)g_{-1}) \bigsqcup (K^+_1\dot+(1,1)g_{-1})
 \bigsqcup(K^+_1\dot+(1,0)g_{-1}\dot+(1,1)g_0)
 $$
 $$
 \bigsqcup(K^+_1\dot+(0,1)g_{-1}\dot+(1,1)g_0).
 $$
 We can consider additive group $K^+$ as product $\mathfrak G\times \mathfrak
 G$ of Cantor groups.
 In this case $\hat\varphi$  and $\varphi$ may be defined  on the product
 $\mathfrak G_{1}^\bot \times \mathfrak G_{1}^\bot$  and
$\mathfrak G_{-1} \times \mathfrak G_{-1}$ respectively
 by the tables

\begin{picture}(90,64)
   \put(0,0){\line(1,0){60}}
   \qbezier(0,0)(15,-4)(30,0)
    \qbezier(0,0)(30,-6)(60,0)
    \qbezier(60,0)(63,8)(60,15)
     \put(11,-8){$\mathfrak G_{0}^\bot$}
     \put(27,-10){$\mathfrak G_{1}^\bot$}
     \put(63,7){$\mathfrak G_{-1}^\bot$}
   \put(0,15){\line(1,0){60}}
   \put(0,30){\line(1,0){60}}
   \put(0,45){\line(1,0){60}}
   \put(0,60){\line(1,0){60}}

   \put(0,0){\line(0,1){60}}
   \put(15,0){\line(0,1){60}}
   \put(30,0){\line(0,1){60}}
   \put(45,0){\line(0,1){60}}
   \put(60,0){\line(0,1){60}}

   \put(5,5){$1$}
   \put(20,20){$1$}
   \put(35,50){$1$}
   \put(50,35){$1$}

   \put(20,5){$0$}  \put(35,5){$0$} \put(50,5){$0$}
   \put(5,20){$0$}  \put(35,20){$0$} \put(50,20){$0$}
   \put(5,35){$0$}  \put(20,35){$0$} \put(35,35){$0$}
   \put(5,50){$0$}  \put(20,50){$0$} \put(50,50){$0$}
  \end{picture}
  \begin{picture}(60,64)
   \put(0,0){\line(1,0){60}}
   \qbezier(0,0)(15,-4)(30,0)
    \qbezier(0,0)(30,-6)(60,0)
    \qbezier(60,0)(63,8)(60,15)
     \put(11,-8){$\mathfrak G_{0}$}
     \put(27,-10){$\mathfrak G_{-1}$}
     \put(63,7){$\mathfrak G_{1}$}
   \put(0,15){\line(1,0){60}}
   \put(0,30){\line(1,0){60}}
   \put(0,45){\line(1,0){60}}
   \put(0,60){\line(1,0){60}}

   \put(0,0){\line(0,1){60}}
   \put(15,0){\line(0,1){60}}
   \put(30,0){\line(0,1){60}}
   \put(45,0){\line(0,1){60}}
   \put(60,0){\line(0,1){60}}

   \put(5,5){$1$}
   \put(20,20){$1$}
   \put(35,50){$1$}
   \put(50,35){$1$}

   \put(20,5){$0$}  \put(35,5){$0$} \put(50,5){$0$}
   \put(5,20){$0$}  \put(35,20){$0$} \put(50,20){$0$}
   \put(5,35){$0$}  \put(20,35){$0$} \put(35,35){$0$}
   \put(5,50){$0$}  \put(20,50){$0$} \put(50,50){$0$}
  \end{picture}

  \vspace*{10mm}

\hspace*{10mm} Figure 6. $\hat\varphi$ \hspace*{40mm} Figure 7.
$\varphi$

 \vspace*{5mm}

 Since ${\rm supp}\hat\varphi \neq
(K_{0}^+)^\bot$ and ${\rm supp}\varphi \neq  (K_{0}^+)$, it
follows that $\varphi$ generates non-Haar MRA. From this example
we see that MRA on local field gives an effective method to
construct multidimensional step wavelets.

\vspace*{5mm}

{\bf Acknowledgements}\\
  The first author was  supported by the  state program
of Russian Ministry of Education and Science (project
1.1520.2014K). The second author was supported by RFBR (grant
13-01-00102).

 \end{document}